\documentclass{degm}

\usepackage{amssymb}
\usepackage{amsmath}
\usepackage{latexsym}
\usepackage{amsthm}
\usepackage{cite}
\usepackage{eucal}

\makeatletter
\def\thmhead@plain#1#2#3{%
  \thmname{#1}\thmnumber{\@ifnotempty{#1}{ }#2}%
  \thmnote{ \the\thm@notefont(#3)}}
\let\thmhead\thmhead@plain
\def\swappedhead#1#2#3{%
  \thmnumber{#2}\thmname{\@ifnotempty{#2}{. }#1}%
  \thmnote{ \the\thm@notefont(#3)}}
\makeatother

\theoremstyle{definition} %%% for statements in roman typeface

\newtheorem{definition}{Definition}[section]
\newtheorem{remark}[definition]{Remark}
\newtheorem{example}[definition]{Example}

\newtheorem{question}[definition]{Question}

\theoremstyle{plain}      %%% for statements in italic typeface

\newtheorem{proposition}[definition]{Proposition}
\newtheorem{theorem}[definition]{Theorem}
\newtheorem{corollary}[definition]{Corollary}
\newtheorem{lemma}[definition]{Lemma}

\title{On the Representation Theory of Deformation Quantization}  

\author{\textbf{Stefan Waldmann\thanks{email:
      Stefan.Waldmann@physik.uni-freiburg.de}}  
  \\[0.5cm]
  Fakult{\"a}t f{\"u}r Physik\\ 
  Albert-Ludwigs-Universit{\"a}t Freiburg\\
  Hermann Herder Stra{\ss}e 3\\
  D 79104 Freiburg\\
  Germany}

\date{June 2001}

\newcommand{\im}         {\mathrm{i}}
\newcommand{\eu}         {\mathrm{e}}
\newcommand{\cc}[1]      {\overline{{#1}}}
\newcommand{\id}         {\mathsf{id}}
\newcommand{\tr}         {\mathsf{tr}}
\newcommand{\End}        {\mathop{{\mathsf{End}}}}
\newcommand{\ring}[1]    {\mathsf{{#1}}}
\newcommand{\SP} [1]     {\left\langle{{#1}}\right\rangle}
\newcommand{\Unit}       {\mathsf{1}}
\newcommand{\rep}        {{}^*\!\mbox{-}\mathsf{rep}}
\newcommand{\Rep}        {{}^*\!\mbox{-}\mathsf{Rep}}
\newcommand{\supp}       {\mathop{{\mathrm{supp}}}}

\newcommand{\Exp}        {\operatorname{\mathrm Exp}}

\newcommand{\qe}         {\boldsymbol{e}}
\newcommand{\qs}         {\boldsymbol{s}}
\newcommand{\qpsi}       {\boldsymbol{\Psi}}
\newcommand{\qphi}       {\boldsymbol{\Phi}}
\newcommand{\qh}         {\boldsymbol{h}}

\newcommand{\qV}         {\boldsymbol{V}}
\newcommand{\qA}         {\boldsymbol{A}}
\newcommand{\qB}         {\boldsymbol{B}}
\newcommand{\qt}         {\boldsymbol{t}}

\newcommand{\CH}         {\check{\mathrm{H}}}
\newcommand{\Pic}        {\mathsf{Pic}}
\newcommand{\HdR}        {\mathrm{H}_{\scriptscriptstyle \mathrm{dR}}}
\newcommand{\Def}        {\mathsf{Def}}

\begin{document}

\maketitle

\begin{abstract}
    In this contribution to the proceedings of the
    $68^{\textrm{\`eme}}$ Rencontre entre Physiciens Th\'eoriciens et
    Math\'ematiciens on Deformation Quantization I shall report on
    some recent joint work with Henrique Bursztyn on the
    representation theory of $^*$-algebras arising from deformation
    quantization as presented in my talk.
\end{abstract}

\section{Introduction}
\label{sec:intro}
The starting point for this investigation was the wish to develop a
representation theory for the star product algebras as they arise in
deformation quantization \cite{bayen.et.al:1978}, see
\cite{gutt:2000,weinstein:1994} for recent reviews as well as Daniel
Sternheimer's contribution in this proceedings. In deformation
quantization the main emphasis lies on the construction of the 
\emph{observable algebra} out of the classical data. The
observables are understood to be the \emph{primary object} as opposed
to other quantization schemes based on particular operator
representations on Hilbert spaces. More specifically, one deforms the
classical Poisson algebra of observables, modeled by smooth
complex-valued functions $C^\infty(M)$ on a classical phase space $M$,
into the direction of the Poisson bracket determined by the symplectic
or Poisson structure of $M$. This new non-commutative product is
understood to be the multiplication law for the quantum
observables. Usually, formal associative deformations in the sense of
Gerstenhaber \cite{gerstenhaber:1964,gerstenhaber.schack:1988} are
considered, yielding the notion of (formal) \emph{star products}. The
$^*$-involution of star products as well as the order structure of
real formal power series allow for an \emph{intrinsic} definition of
\emph{states} as positive linear functionals of the observable
algebra. It turns out that many of the well-known `quantizations' can
be recovered as GNS representations induced by certain particular
states 
\cite{bordemann.neumaier.waldmann:1999,bordemann.waldmann:1998,bordemann.neumaier.waldmann:1998,bordemann.neumaier.pflaum.waldmann:1998:pre,bordemann.roemer.waldmann:1998,waldmann:2000}.
The reader will certainly notice that this line of arguments reminds
very much the Haag-Kastler formulation of quantum field theory
\cite{haag:1993}.

Let us now discuss the motivation and conceptual advantages for such an
approach to quantization: First the states are known automatically
(and even in some functorial way) once the observable algebra is
known. Secondly, and even more important, one can consider
\emph{different} representations of one observable
algebra. This corresponds to the same system in
\emph{different} physical situations. Here one has various examples
like thermodynamical equilibrium at different temperatures, modeled
e.g.\  by KMS states
\cite{basart.et.al:1984,bordemann.roemer.waldmann:1998} in order to
study phase transitions, as well as
super selection rules and topological effects like the Aharonov-Bohm
effect \cite{bordemann.neumaier.pflaum.waldmann:1998:pre}. 
A further example is the re-interpretation of Dirac's quantization
condition for magnetic charges of magnetic monopoles as a condition
for Morita equivalence of the quantization with and without magnetic
field. In particular, the whole categories of $^*$-representations are
equivalent when Dirac's condition is fulfilled
\cite{bursztyn.waldmann:2001:pre}.

These examples show the need for a good understanding of what the
representation theory of deformation quantization looks like.
Thus we shall give here an overview over the recent developments and
discuss various examples. We shall focus on the ideas and motivations
rather than on the detailed proofs as they can be found
elsewhere 
\cite{bursztyn.waldmann:2001a,bursztyn:2001:pre,bursztyn.waldmann:2001:pre,bursztyn.waldmann:2000b,bursztyn.waldmann:2000a:pre,bordemann.waldmann:1998,bursztyn.waldmann:2000a}. 
We also give a quite extensive (though certainly incomplete) guide to
the literature.

The article is organized as follows: In Sect.~\ref{sec:defquant} we
recall the basic notions of deformation quantization which provides the
main example. In Sect.~\ref{sec:rings} the concept of $^*$-algebras
over ordered rings and their $^*$-representations is discussed and
Sect.~\ref{sec:gns} is devoted to the fundamental GNS construction of
$^*$-representations. Then in Sect.~\ref{sec:defalg} we point out some
questions and results on the deformation theory of
$^*$-algebras. Sect.~\ref{sec:rieffel} contains a brief description of
Rieffel induction and in Sect.~\ref{sec:morita} we establish the
notion of Morita equivalence of $^*$-algebras. In
Sect.~\ref{sec:module} we discuss the deformation of projective
modules as they provide the main example of Morita equivalence
bimodules. This will be applied to Hermitian vector bundles in 
Sect.~\ref{sec:vector} where we obtain local expressions for the
deformations. Sect.~\ref{sec:ltog} will give the way back, namely the
construction of global deformations out of local data. Finally, we
compute the relative characteristic class of Morita equivalent star
products on symplectic manifolds in Sect.~\ref{sec:mestar}.

\noindent
\textbf{Acknowledgments:}
First of all it is a pleasure for me to thank the organizers and in
particular Gilles Halbout for the invitation to present these results
at the $68^{\textrm{\`eme}}$ Rencontre entre Physiciens Th\'eoriciens
et Math\'ematiciens on `Deformation Quantization' in
Strasbourg. I also would like to thank Martin Bordemann, Henrique
Bursztyn, Murray Gerstenhaber, Simone Gutt for valuable discussions
and remarks. Moreover, I'm very grateful to Branislav Jur\v{c}o and
Peter Schupp for a discussion clarifying the appearance of the
diffeomorphism in Theorem~\ref{theorem:mestars}.
\section{Deformation Quantization as Example}
\label{sec:defquant}
Let us briefly recall the notion of deformation quantization and star
products \cite{bayen.et.al:1978} in order to fix our notation. 
\begin{definition}
    \label{definition:starproduct}
    Let $(M, \pi)$ be a Poisson manifold. A \emph{star product}
    $\star$ on $(M, \pi)$ is a formal associative
    $\mathbb{C}[[\lambda]]$-bilinear product for
    $C^\infty(M)[[\lambda]]$, 
    \begin{equation}
        \label{eq:starproduct}
        f \star g = \sum_{r=0}^\infty \lambda^r C_r (f, g),
    \end{equation}
    with bidifferential operators $C_r$ such that $C_0(f,g) = fg$ and
    $C_1(f,g) - C_1(g,f) = \im \{f, g\}$ as well as 
    $1 \star f = f = f \star 1$. The star product $\star$ is called
    \emph{Hermitian} if $\cc{f \star g} = \cc{g} \star \cc{f}$, where
    $\cc{\lambda} = \lambda$. Two star products $\star$, $\star'$ are
    called \emph{equivalent} if there exists a formal series of
    differential operators $T = \id + \sum_{r=1}^\infty \lambda^r T_r$
    such that 
    \begin{equation}
        \label{eq:equiv}
        T(f \star' g) = Tf \star Tg.
    \end{equation}
    Hermitian star products are \emph{$^*$-equivalent} if they are
    equivalent with a $T$ such that $T(\cc{f}) = \cc{T(f)}$.
\end{definition}
With the interpretation $\lambda \leftrightarrow \hbar$ one identifies
$C^\infty(M)[[\lambda]]$ with the algebra of observables of the
quantum system corresponding to $(M, \pi)$. However, there are
applications of star products beyond quantization of phase spaces
coming from certain limits of string
theory\cite{schomerus:1999,jurco.schupp.wess:2001a:pre,jurco.schupp.wess:2001b:pre,jurco.schupp.wess:2000,connes.douglas.schwarz:1998,waldmann:2001a:pre}.
Here they give rise to non-commutative versions of space-time very much in
the spirit of Connes' non-commutative geometry \cite{connes:1994}.

The existence of star products was first established in the symplectic
case
\cite{dewilde.lecomte:1983b,fedosov:1994a,omori.maeda.yoshioka:1991}
and recently also in the general Poisson case
\cite{kontsevich:1997:pre,cattaneo.felder:2000}. Moreover, the
classification up to equivalence is also well-understood
\cite{nest.tsygan:1995a,bertelson.cahen.gutt:1997,gutt.rawnsley:1999,kontsevich:1997:pre}.

Having a quantized observable algebra 
$\boldsymbol{\mathcal{A}} = (C^\infty(M)[[\lambda]],\star)$ we want to
understand the states as well as the representations of
$\boldsymbol{\mathcal{A}}$ on Hilbert spaces. At first glance
this looks almost hopeless due to the formal power series. There are
no interesting $\mathbb{C}$-linear functionals
$\omega: C^\infty (M)[\lambda]] \to \mathbb{C}$ which could serve as
states nor reasonable representations on Hilbert spaces. One
immediately is confronted with the delicate question of convergence.

This changes drastically as soon as one takes the formal series
seriously and considers $\mathbb{C}[[\lambda]]$-linear functionals
\begin{equation}
    \label{eq:states}
    \omega: C^\infty(M)[[\lambda]] \to \mathbb{C}[[\lambda]]
\end{equation}
instead. Now we observe that $\mathbb{R}[[\lambda]]$ is an
\emph{ordered} ring (see Remark~\ref{remark:ordered}) whence the
positivity condition 
\begin{equation}
    \label{eq:posstate}
    \omega( \cc{f} \star f) \ge 0
\end{equation}
makes sense as soon as $\star$ is a Hermitian star product. This way
we indeed obtain a notion of positive linear functionals and thus for
states, now completely within the framework of formal
deformations. Moreover, as many examples show, this notion gives
physically reasonable results and allows for a mathematically
non-trivial theory of $^*$-representations which we shall discuss in
the sequel.
\section{$^*$-Algebras over Ordered Rings}
\label{sec:rings}
Having deformation quantization as the major `example' we shall now
focus on the general algebraic structures in order to develop a more
axiomatic approach to a representation theory for star products.

Let $\ring{R}$ be an ordered ring, i.e.\  a commutative associative
unital ring with a distinguished subset $\ring{P} \subset \ring{R}$ of
\emph{positive elements} such that $\ring{R}$ is the disjoint union
$\ring{R} = -\ring{P} \cup \{0\} \cup \ring{P}$ and
$\ring{P} \cdot \ring{P} \subseteq \ring{P}$ as well as 
$\ring{P} + \ring{P} \subseteq \ring{P}$. As usual we shall write 
$a < b$ if $b - a \in \ring{P}$ etc. Moreover, we set
$\ring{C} = \ring{R}(\im)$ with $\im^2 = -1$. As consequences we note
that $\ring{R}$ as well as $\ring{C}$ have characteristic zero,
i.e. $\mathbb{Z} \subseteq \ring{R} \subset \ring{C}$. Moreover,
$\ring{R}$ and $\ring{C}$ have no zero divisors and the quotient field
$\hat{\ring{R}}$ of $\ring{R}$ is canonically an ordered field. The
complex conjugation $z \mapsto \cc{z}$ is an involutive ring
automorphism of $\ring{C}$ with $\cc{z}z > 0$ if and only if 
$z \ne0$. Here we use the fact that $\ring{R}$ is embedded 
into $\ring{C}$ in the usual way.
\begin{remark}
    \label{remark:ordered}
    If $\ring{R}$ is ordered then $\ring{R}[[\lambda]]$ is
    ordered in a canonical way, too, by the definition 
    $\sum_{r=r_0}^\infty \lambda^r a_r > 0$ if
    $a_{r_0} > 0$. This shows that the concept of ordered rings fits
    naturally to formal power series and thus to Gerstenhaber's
    deformation theory \cite{gerstenhaber.schack:1988}. 
    Then $\ring{R}[[\lambda]]$ will be \emph{non-Archimedian} since
    e.g. $n\lambda < 1$ for all $n \in \mathbb{Z}$. The interpretation
    in formal deformation theory is that the deformation parameter
    $\lambda$ is `very small' compared to the other numbers in
    $\ring{R}$.
\end{remark}

In the following we shall use a fixed choice $\ring{R}$ and $\ring{C}$
as replacement for $\mathbb{R}$ and $\mathbb{C}$. The order structure
allows for the following definition:
\begin{definition}
    \label{definition:prehilbert}
    A $\ring{C}$-module $\mathfrak{H}$ is called 
    \emph{pre-Hilbert space} if it is equipped with an inner product
    $\SP{\cdot,\cdot}: 
    \mathfrak{H} \times \mathfrak{H} \to \ring{C}$ such that
    $\SP{\phi, z\psi + w\chi} = z\SP{\phi,\psi} + w \SP{\phi,\chi}$
    and $\SP{\phi, \psi} = \cc{\SP{\psi,\phi}}$ for $z,w \in \ring{C}$
    and $\phi,\psi,\chi \in \mathfrak{H}$ as well as 
    $\SP{\phi,\phi} > 0$ for $\phi \ne 0$.
\end{definition}
The simplest example is clearly $\mathfrak{H} = \ring{C}^n$, where the
inner product is defined as usual by
\begin{equation}
    \label{eq:standardinnerproduct}
    \SP{z,w} = \sum_{i=1}^n \cc{z}_i w_i.
\end{equation}

Having a pre-Hilbert space over $\ring{C}$ we need an algebraic
replacement for the `bounded operators'. Since in the ordinary theory
of Hilbert spaces over the complex numbers an operator is bounded if
and only if it is adjointable (Hellinger-Toepliz theorem, see
e.g.~\cite[p.~117]{rudin:1991}), we state the following definition:
\begin{definition}
    \label{definition:adjoints}
    Let $\mathfrak{H}, \mathfrak{K}$ be pre-Hilbert spaces over 
    $\ring{C}$. A $\ring{C}$-linear map 
    $A: \mathfrak{H} \to \mathfrak{K}$ is called \emph{adjointable} if
    there exists a $\ring{C}$-linear map
    $A^*: \mathfrak{K} \to \mathfrak{H}$ such that
    \begin{equation}
        \label{eq:adoint}
        \SP{A \phi, \psi} = \SP{\phi, A^*\psi}
    \end{equation}
    for all $\phi \in \mathfrak{H}$ and $\psi \in \mathfrak{K}$. The
    set of adjointable maps is denoted by
    $\mathfrak{B}(\mathfrak{H},\mathfrak{K})$ and
    $\mathfrak{B}(\mathfrak{H}) :=
    \mathfrak{B}(\mathfrak{H},\mathfrak{H})$. Moreover,
    $U \in \mathfrak{B}(\mathfrak{H}, \mathfrak{K})$ is called
    \emph{isometric} if $\SP{U\phi, U\psi} = \SP{\phi, \psi}$ for all
    $\phi, \psi \in \mathfrak{H}$, and $U$ is called \emph{unitary} if
    it is isometric and bijective.
\end{definition}
\begin{lemma}
    \label{lemma:adjoints}
    Let $\mathfrak{H}, \mathfrak{K}, \mathfrak{L}$ be pre-Hilbert
    spaces over $\ring{C}$ and let 
    $A, B \in \mathfrak{B}(\mathfrak{H},\mathfrak{K})$,
    $C \in \mathfrak{B}(\mathfrak{K},\mathfrak{L})$, and 
    $z, w \in \ring{C}$. 
    Then the adjoint $A^*$ of $A$ is unique and we have the usual
    rules 
    \begin{equation}
        \label{eq:adjoints}
        (zA + wB)^* = \cc{z}A^* + \cc{w}B^*,
        \quad
        (A^*)^* = A,
        \quad
        (CA)^* = A^*C^*,
    \end{equation}
    where the existence of the corresponding adjoints is
    automatic. Moreover, 
    $U \in \mathfrak{B}(\mathfrak{H}, \mathfrak{K})$ is isometric if
    and only if $U^*U = \id_{\mathfrak{H}}$ and unitary if and only if
    $U^* = U^{-1}$.
\end{lemma}
Thus the pre-Hilbert spaces over $\ring{C}$ constitute a category
$\textsf{pre-Hilbert}(\ring{C})$ where the morphisms are given by
adjointable $\ring{C}$-linear maps. Directly related to pre-Hilbert
spaces are the $^*$-algebras over $\ring{C}$.
\begin{definition}
    \label{definition:staralg}
    An associative $\ring{C}$-algebra $\mathcal{A}$ is called
    \emph{$^*$-algebra} over $\ring{C}$ if it is equipped with a
    $^*$-involution, i.e.\  an anti-linear, involutive anti-automorphism
    $A \mapsto A^*$. For two $^*$-algebras $\mathcal{A}$,
    $\mathcal{B}$ a \emph{$^*$-homomorphism}
    $\Phi: \mathcal{A} \to \mathcal{B}$ is a homomorphism of
    associative algebras such that $\Phi(A^*) = \Phi(A)^*$ for all 
    $A \in \mathcal{A}$.    
\end{definition}
\begin{example}[$^*$-Algebras]
    \label{example:staralgs}
    \hfill
    \begin{enumerate}
    \item Primary example and guideline for the following are the
        $C^*$-algebras over $\mathbb{C}$.
    \item The main motivation for generalizing the theory of
        $C^*$-algebras to this more algebraic framework are the star
        products from deformation quantization. Here one also has
        the non-unital versions $(C^\infty_0(M)[[\lambda]], \star)$.
    \item If $\mathfrak{H}$ is a pre-Hilbert space over $\ring{C}$
        then $\mathfrak{B}(\mathfrak{H})$ is a unital $^*$-algebra
        over $\ring{C}$. In particular, if 
        $\mathfrak{H} = \ring{C}^n$ with inner product
        \eqref{eq:standardinnerproduct} we have
        $\mathfrak{B}(\ring{C}^n) \cong M_n (\ring{C})$
        with the usual $^*$-involution.
    \item For a pre-Hilbert space $\mathfrak{H}$ we define the finite
        rank operators $\mathfrak{F}(\mathfrak{H})$ to be the
        $\ring{C}$-linear span of the operators $\Theta_{\phi,\psi}$
        which are defined by 
        $\Theta_{\phi,\psi}(\chi) = \phi \SP{\psi,\chi}$ for
        $\phi,\psi,\chi \in \mathfrak{H}$. Then
        $\mathfrak{F}(\mathfrak{H}) \subseteq
        \mathfrak{B}(\mathfrak{H})$ 
        is a $^*$-ideal and hence a $^*$-algebra itself. It is a
        proper ideal in infinite dimensions while
        $\mathfrak{F}(\ring{C}^n) = \mathfrak{B}(\ring{C}^n)
        \cong M_n (\ring{C})$.
    \item If $\mathcal{A}$ is a $^*$-algebra over $\ring{C}$ then
        $M_n(\mathcal{A})$ with the usual product and $^*$-involution
        is again a $^*$-algebra over $\ring{C}$.
    \item There are various examples beyond the above ones as
        e.g.\  algebras of differential operators or complexified
        universal enveloping algebras
        $U(\mathfrak{g}) \otimes_{\ring{R}} \ring{C}$ of Lie algebras
        $\mathfrak{g}$ over $\ring{R}$ with $^*$-involution induced by
        $x \mapsto -x$ for $x \in \mathfrak{g}$.
    \end{enumerate}
\end{example}
As usual we define \emph{Hermitian}, \emph{isometric}, \emph{unitary},
and \emph{normal elements} of a $^*$-algebra. The order
structure of $\ring{R}$ becomes crucial when we want to discuss the
following notions of positivity:
\begin{definition}
    \label{definition:positivities}
    Let $\mathcal{A}$ be a $^*$-algebra over $\ring{C}$.
    A linear functional 
    $\omega: \mathcal{A} \to \ring{C}$ is called \emph{positive} if
    \begin{equation}
        \label{eq:posfun}
        \omega(A^*A) \ge 0
    \end{equation}
    for all $A \in \mathcal{A}$. If $\mathcal{A}$ is unital
    $\omega$ is called a \emph{state} if in addition 
    $\omega(\Unit) = 1$.
    An algebra element $A \in \mathcal{A}$ is called
    \emph{algebraically positive} if 
    $A = b_1 B_1^*B_1 + \cdots + b_n B_n^*B_n$ with
    $b_i > 0$ and $B_i \in \mathcal{A}$. 
    Moreover, $A = A^* \in \mathcal{A}$ is called \emph{positive} if
    for all positive linear functionals 
    $\omega: \mathcal{A} \to \ring{C}$ one has $\omega(A) \ge 0$.
\end{definition}
The algebraically positive elements are denoted by $\mathcal{A}^{++}$
while the positive ones are denoted by $\mathcal{A}^+$. Clearly
$\mathcal{A}^{++} \subseteq \mathcal{A}^+$ and both subsets are convex
cones in $\mathcal{A}$. A positive linear functional 
$\omega: \mathcal{A} \to \ring{C}$ is real in the sense that
$\omega(A^*B) = \cc{\omega(B^*A)}$, which implies 
$\omega(A^*) = \cc{\omega(A)}$ in case $\mathcal{A}$ is
unital. Moreover, we have the \emph{Cauchy-Schwarz} inequality
\begin{equation}
    \label{eq:CSU}
    \omega(A^*B) \cc{\omega(A^*B)} \le \omega(A^*A) \omega(B^*B)
\end{equation}
for all $A, B \in \mathcal{A}$. Clearly positive linear functionals
can be pulled back under $^*$-homomorphisms whence $^*$-homomorphisms
map (algebraically) positive elements to (algebraically) positive
elements.
\begin{definition}
    \label{definition:starrep}
    Let $\mathcal{A}$ be a $^*$-algebra over $\ring{C}$. A
    \emph{$^*$-representation} $\pi$ of $\mathcal{A}$ on a pre-Hilbert
    space $\mathfrak{H}$ is a $^*$-homomorphism 
    $\pi: \mathcal{A} \to \mathfrak{B}(\mathfrak{H})$. If 
    $(\pi,\mathfrak{H})$ and $(\rho, \mathfrak{K})$ are
    $^*$-representations of $\mathcal{A}$ then a $\ring{C}$-linear map
    $U: \mathfrak{H} \to \mathfrak{K}$ is called \emph{intertwiner} if
    \begin{equation}
        \label{eq:intertwiner}
        U \pi(A) = \rho(A) U
    \end{equation}
    for $A \in \mathcal{A}$. The category of $^*$-representations of
    $\mathcal{A}$ with isometric intertwiners as morphisms is denoted
    by $\rep(\mathcal{A})$. The category of strongly non-degenerate
    $^*$-representations 
    (i.e. $\pi(\mathcal{A})\mathfrak{H} = \mathfrak{H}$) of
    $\mathcal{A}$ is denoted by $\Rep(\mathcal{A})$. 
\end{definition}
If $\mathcal{A}$ is unital, strong non-degeneracy is equivalent
to $\pi(\Unit) = \id$. We shall be mainly interested in the unital
case in the following.

However, there are two notions which generalize
the unital case in a useful way. An algebra is called
\emph{idempotent} if any element can be obtained as sum of products 
$A = B_1C_1 + \cdots + B_nC_n$. Even more useful is the notion of an
\emph{approximate identity}: Let $I$ be a directed set and let
$\mathcal{A}_\alpha \subseteq \mathcal{A}$ be sub-algebras and 
$E_\alpha \in \mathcal{A}$ such that 
$\mathcal{A}_\alpha \subseteq \mathcal{A}_\beta$ for
$\alpha \le \beta$, $\mathcal{A} = \cup_\alpha \mathcal{A}_\alpha$,
$E_\alpha E_\beta = E_\alpha = E^*_\alpha = E_\beta E_\alpha$ for all
$\alpha < \beta$ and $E_\alpha A = A = A E_\alpha$ for all 
$A \in \mathcal{A}_\alpha$. Then
$\{\mathcal{A}_\alpha,E_\alpha\}_{\alpha \in I}$ is called an
approximate identity of $\mathcal{A}$. Note that we do not require
$E_\alpha^2 = E_\alpha$ nor $E_\alpha \in \mathcal{A}_\alpha$. As an
example on can consider $C^\infty_0(M)$ for $M$ non-compact: choose
open subsets $O_n \subseteq M$ with compact closure such that
$O_n^{\mathrm{cl}} \subseteq O_{n+1}$ and $M = \cup_n O_n$. Moreover,
choose $\chi_n = \cc{\chi}_n \in C^\infty_0 (O_{n+1})$ such that
$\chi_n \big|_{O_n^{\mathrm{cl}}} = 1$. Then 
$\{C^\infty_0(O_n), \chi_n\}_{n \in \mathbb{N}}$ is a countable
approximate identity for $C^\infty_0(M)$.

Most of the results for unital $^*$-algebras can be
extended to the non-unital case provided the $^*$-algebra has an
approximate identity. Nevertheless, we shall not emphasize this aspect
too much but refer to
\cite{bursztyn.waldmann:2001a,bursztyn.waldmann:2000a:pre} for the
corresponding statements.
\section{The GNS Construction}
\label{sec:gns}
In general, there may be no interesting $^*$-representations of
$\mathcal{A}$ at all. However, the order structure of $\ring{R}$
allows for the GNS construction of $^*$-representations out of
positive linear functionals.

We shall now describe this construction. Let 
$\omega: \mathcal{A} \to \ring{C}$ be positive and define
\begin{equation}
    \label{eq:Gelfand}
    \mathcal{J}_\omega 
    = \{A \in \mathcal{A} \; | \; \omega(A^*A) = 0\}.
\end{equation}
The \eqref{eq:CSU} implies that $\mathcal{J}_\omega$ is a left ideal
of $\mathcal{A}$, called the \emph{Gel'fand ideal} of $\omega$. Hence
\begin{equation}
    \label{eq:GNSpreHilbert}
    \mathfrak{H}_\omega = \mathcal{A} \big/ \mathcal{J}_\omega
\end{equation}
becomes a $\mathcal{A}$-left module with the usual left action denoted
by $\pi_\omega (A)\psi_B  = \psi_{AB}$, where 
$\psi_B \in \mathfrak{H}_\omega$ is the class of $B \in
\mathcal{A}$. Now $\mathfrak{H}_\omega$ becomes a pre-Hilbert
space by
\begin{equation}
    \label{eq:GNSinnerprod}
    \SP{\psi_B,\psi_C}_\omega = \omega(B^*C).
\end{equation}
Finally, it is verified easily that 
$\pi_\omega (A) \in \mathfrak{B}(\mathfrak{H}_\omega)$ and
$\pi_\omega(A^*) = \pi_\omega(A)^*$ whence $\pi_\omega$ is a
$^*$-representation of $\mathcal{A}$ on $\mathfrak{H}_\omega$, the
so-called \emph{GNS representation} induced by $\omega$. In case
$\mathcal{A}$ is unital one even has a \emph{cyclic vector}, namely
$\psi_{\Unit}$, whence $\pi_\omega$ is a 
\emph{cyclic $^*$-representation} and hence strongly
non-degenerate. In this case
\begin{equation}
    \label{eq:omegafrompiomega}
    \omega(A) = \SP{\psi_{\Unit}, \pi_\omega(A) \psi_{\Unit}}.
\end{equation}
More generally, if $\mathcal{A}$ is idempotent then $\pi_\omega$ is
still strongly non-degenerate. 
\begin{remark}
    In the case of $C^*$-algebras this construction is of course a
    classical result, see e.g.\  the textbook
    \cite{landsman:1998}. The surprising 
    observation is that it can be generalized to the above purely
    algebraic framework and still gives non-trivial examples: Many
    `operator representations' by 
    differential operators turn out to be particular GNS
    representations of certain star product algebras. We shall not
    discuss all the examples here but refer to the literature: See
    \cite{bordemann.neumaier.waldmann:1998,bordemann.neumaier.waldmann:1999,bordemann.neumaier.pflaum.waldmann:1998:pre} 
    for Schr\"odinger type representations on cotangent bundles,
    \cite{bordemann.waldmann:1998} for Bargmann-Fock type
    representations for Fedosov star products of Wick type on K\"ahler
    manifolds \cite{karabegov:2000,karabegov:1996,bordemann.waldmann:1997a},
    and \cite{waldmann:2000}
    for the left regular representation induced by KMS functionals
    \cite{basart.et.al:1984,bordemann.roemer.waldmann:1998}. 
\end{remark}

With the GNS construction in mind it is not very surprising that the
question of existence of `interesting' (e.g.\  faithful)
$^*$-representations is intimately linked to the existence of positive
linear functionals.
\begin{definition}
    \label{definition:suff}
    A $^*$-algebra $\mathcal{A}$ over $\ring{C}$ has 
    \emph{sufficiently many positive linear functionals} if for every
    non-zero Hermitian element $H \in \mathcal{A}$ there exists a
    positive linear functional $\omega: \mathcal{A} \to \ring{C}$ such
    that $\omega(H) \ne 0$.
\end{definition}
There are simple examples of $^*$-algebras \emph{without} sufficiently
many positive linear functionals, see
e.g.~\cite[Sect.~2]{bursztyn.waldmann:2000a:pre}.
\begin{definition}
    \label{definition:closedideals}
    Let $\mathcal{A}$ be a $^*$-algebra over $\ring{C}$ and
    $\mathcal{J} \subseteq \mathcal{A}$ a $^*$-ideal. Then
    $\mathcal{J}$ is called \emph{closed} if it is the kernel of a
    $^*$-representation. Moreover, the \emph{minimal ideal} of
    $\mathcal{A}$ is defined by
    \begin{equation}
        \label{eq:minimalideal}
        \mathcal{J}_{\mathrm{min}}(\mathcal{A}) =
        \bigcap_{\mathcal{J}\:\textrm{closed}} \mathcal{J}.
    \end{equation}
\end{definition}
As expected from the $C^*$-algebra case the closed $^*$-ideals form a
lattice with minimal element
$\mathcal{J}_{\mathrm{min}}(\mathcal{A})$. The following theorem
relates the existence of positive linear functionals to the `size' of
$\mathcal{J}_{\mathrm{min}}(\mathcal{A})$ and the existence of
$^*$-representations.
\begin{theorem}
    \label{theorem:suff}
    Let $\mathcal{A}$ be a unital $^*$-algebra over $\ring{C}$. Then
    the following statements are equivalent:
    \begin{enumerate}
    \item $\mathcal{A}$ has sufficiently many positive linear
        functionals.
    \item $\mathcal{A}$ has a faithful $^*$-representation.
    \item $\mathcal{J}_{\mathrm{min}}(\mathcal{A}) = \{0\}$.
    \end{enumerate}
    In general one has the characterization
    \begin{equation}
        \label{eq:minimal}
        \mathcal{J}_{\mathrm{min}}(\mathcal{A})
        = \bigcap_\omega \ker \pi_\omega
        = \bigcap_\omega \mathcal{J}_\omega
        = \bigcap_\omega \ker \omega,
    \end{equation}
    where the intersection runs over all positive functionals of
    $\mathcal{A}$. Moreover, if $A^*A = 0$, or if $A$ is normal and
    nilpotent, or if $zA = 0$ for $0 \ne z \in \ring{C}$ then 
    $A \in \mathcal{J}_{\mathrm{min}}(\mathcal{A})$.
\end{theorem}
The proof relies on the GNS construction, see
\cite[Prop.~2.8]{bursztyn.waldmann:2001a} and 
\cite[Thm.~4.7]{bursztyn.waldmann:2000a:pre}. In particular, if there
exists a faithful $^*$-representation then one can already obtain a
faithful $^*$-representation by direct sums of GNS
representations. Moreover, since all the unpleasant and
`non-$C^*$-algebra like' elements of $\mathcal{A}$ are contained in
$\mathcal{J}_{\mathrm{min}}(\mathcal{A})$ we can pass to the quotient
\begin{equation}
    \label{eq:AmodJmin}
    \varrho: 
    \mathcal{A} \to 
    \mathcal{A} \; \big/ \; \mathcal{J}_{\mathrm{min}}(\mathcal{A}).
\end{equation}
This procedure is clearly functorial and can be seen as an algebraic
analog of the construction of the \emph{$C^*$-enveloping algebra} of
a Banach $^*$-algebra, see e.g.~\cite[Sect.~II.7]{dixmier:1977b}.
In particular, the pull-back $\varrho^*$ of $^*$-representations of 
$\mathcal{A} \big/ \mathcal{J}_{\mathrm{min}}(\mathcal{A})$ to
$\mathcal{A}$ implements an equivalence of the categories
$\Rep(\mathcal{A})$ and
$\Rep(\mathcal{A} \big/ \mathcal{J}_{\mathrm{min}}(\mathcal{A}))$. 
Thus elements in $\mathcal{J}_{\mathrm{min}}(\mathcal{A})$ are
`invisible' in $^*$-representations of $\mathcal{A}$, see
\cite[Prop.~5.6]{bursztyn.waldmann:2000a:pre}.
\section{Deformation of $^*$-Algebras}
\label{sec:defalg}
Now we shall consider formal deformations of $^*$-algebras. Since
$\ring{R}[[\lambda]]$ is again ordered one stays in the same framework
of $^*$-algebras over ordered rings.
\begin{definition}
    \label{definition:hermdef}Let $\mathcal{A}$ be a $^*$-algebra over
    $\ring{C}$ and let 
    $\boldsymbol{\mathcal{A}} = (\mathcal{A}[[\lambda]], \star)$ be a
    formal associative deformation of $\mathcal{A}$, viewed as an
    algebra over $\ring{C}[[\lambda]]$. Then
    $\boldsymbol{\mathcal{A}}$ is called a 
    \emph{Hermitian deformation} if 
    \begin{equation}
        \label{eq:hermdef}
        (A \star B)^* = B^* \star A^*.
    \end{equation}
\end{definition}
In principle one can also deform the $^*$-involution of $\mathcal{A}$
in order to obtain a $^*$-involution for $\star$ but we shall only
consider the case of Hermitian deformations here, see also
\cite[Sect.~8]{bursztyn.waldmann:2001a}. Hermitian star products are
of course our main example. For later use we shall recall the
following lemma \cite[Lem.~2.1]{bursztyn.waldmann:2000b}.
\begin{lemma}
    \label{lemma:deformsquares}
    Let $B_0 \in M_n(\mathcal{A})$ be invertible and let 
    $A = B_0^*B_0 + \sum_{r=1}^\infty \lambda^r A_r \in
    M_n(\mathcal{A})[[\lambda]]$. Let $\star$ be a Hermitian
    deformation of $\mathcal{A}$. Then there exists 
    $B_1, \ldots \in M_n(\mathcal{A})$ such that for 
    $B = \sum_{r=0}^\infty \lambda^r B_r$ one has 
    $A = B^* \star B$.
\end{lemma}

In order to compare $^*$-representations of $\mathcal{A}$ and
$\boldsymbol{\mathcal{A}}$ we first need the notion of the
\emph{classical limit} for pre-Hilbert spaces, see 
\cite[Lem.~8.2 and 8.3]{bursztyn.waldmann:2001a}.
\begin{lemma}
    \label{lemma:classlimhilbert}
    Let $\boldsymbol{\mathfrak{H}}$ be a pre-Hilbert space over
    $\ring{C}[[\lambda]]$. Then
    \begin{equation}
        \label{eq:classnull}
        \boldsymbol{\mathfrak{H}}_0 
        = \{\boldsymbol{\phi} \in \boldsymbol{\mathfrak{H}} \; | \; 
        \SP{\boldsymbol{\phi}, \boldsymbol{\phi}}|_{\lambda = 0}
        = 0 \}
    \end{equation}
    is a $\ring{C}[[\lambda]]$ sub-module and 
    $\mathfrak{H} := \mathfrak{C}\boldsymbol{\mathfrak{H}} 
    = \boldsymbol{\mathfrak{H}} \big/ \boldsymbol{\mathfrak{H}}_0$ 
    becomes a pre-Hilbert space over $\ring{C}$ by setting
    $\SP{\mathfrak{C}\boldsymbol{\phi}, \mathfrak{C}\boldsymbol{\psi}}
    = \SP{\boldsymbol{\phi}, \boldsymbol{\psi}}|_{\lambda = 0}$.
    Moreover, if 
    $\boldsymbol{A} \in
    \mathfrak{B}(\boldsymbol{\mathfrak{H}},\boldsymbol{\mathfrak{K}})$
    then 
    $\mathfrak{C}\boldsymbol{A}:\mathfrak{C}\boldsymbol{\mathfrak{H}}
    \to \mathfrak{C}\boldsymbol{\mathfrak{K}} $ 
    defined by
    \begin{equation}
        \label{eq:classlimmap}
        \mathfrak{C}\boldsymbol{A} (\mathfrak{C}\boldsymbol{\phi})
        = \mathfrak{C} (\boldsymbol{A}\boldsymbol{\phi})
    \end{equation}
    is again adjointable 
    $\mathfrak{C}\boldsymbol{A} \in
    \mathfrak{B}(\mathfrak{C}\boldsymbol{\mathfrak{H}},
    \mathfrak{C}\boldsymbol{\mathfrak{K}})$.
\end{lemma}
It follows immediately that $\mathfrak{C}$ defines a functor, the
so-called \emph{classical limit}
\begin{equation}
    \label{eq:classlimfunct}
    \mathfrak{C}: 
    \textsf{pre-Hilbert}(\ring{C}[[\lambda]])
    \longrightarrow \textsf{pre-Hilbert}(\ring{C}).
\end{equation}
\begin{example}
    \label{example:classlimhilbert}
    Let $\mathfrak{H}$ be a pre-Hilbert space over $\ring{C}$. Then
    $\boldsymbol{\mathfrak{H}} = \mathfrak{H}[[\lambda]]$ is a
    pre-Hilbert space over $\ring{C}[[\lambda]]$ by $\lambda$-linear
    extension of the inner product such that canonically
    $\mathfrak{C}\boldsymbol{\mathfrak{H}} \cong
    \mathfrak{H}$. Moreover, 
    $\mathfrak{B}(\boldsymbol{\mathfrak{H}}) =
    \mathfrak{B}(\mathfrak{H})[[\lambda]]$. 
\end{example}
This example shows that the classical limit functor $\mathfrak{C}$ is
essentially surjective. On the other hand there are pre-Hilbert spaces
$\boldsymbol{\mathfrak{H}}$ over $\ring{C}[[\lambda]]$ which are not
of the form $\mathfrak{H}[[\lambda]]$ as in
Example~\ref{example:classlimhilbert}, even if
$\boldsymbol{\mathfrak{H}} \cong \mathfrak{H}[[\lambda]]$ as
$\ring{C}[[\lambda]]$-modules, since the inner products may depend
on $\lambda$ in a non-trivial way.

One also can compute the classical limit of $^*$-representations of
Hermitian deformations in a functorial way, see
\cite[Prop.~8.5]{bursztyn.waldmann:2001a}.
\begin{lemma}
    Let $\boldsymbol{\mathcal{A}}$ be a Hermitian deformation of a
    $^*$-algebra $\mathcal{A}$ over $\ring{C}$ and let
    $\boldsymbol{\pi}: \boldsymbol{\mathcal{A}} \to
    \mathfrak{B}(\boldsymbol{\mathfrak{H}})$ 
    be a $^*$-representation. Then
    $(\mathfrak{C}\boldsymbol{\pi})(A) 
    = \mathfrak{C}(\boldsymbol{\pi}(A))$
    defines a $^*$-representation of $\mathcal{A}$ on 
    $\mathfrak{H} = \mathfrak{C}\boldsymbol{\mathfrak{H}}$ and
    \begin{equation}
        \label{eq:classlimfunctrep}
        \mathfrak{C}: \rep(\boldsymbol{\mathcal{A}}) 
        \longrightarrow \rep(\mathcal{A})
    \end{equation}
    is functorial such that
    $\mathfrak{C}(\Rep(\boldsymbol{\mathcal{A}})) 
    \subseteq \Rep(\mathcal{A})$.    
\end{lemma}
\begin{definition}
    Let $(\pi, \mathfrak{H})$ be a $^*$-representation of
    $\mathcal{A}$ and $\boldsymbol{\mathcal{A}}$ a Hermitian
    deformation of $\mathcal{A}$. A \emph{deformation} of
    $(\pi, \mathfrak{H})$ with respect to $\boldsymbol{\mathcal{A}}$
    is a $^*$-representation
    $(\boldsymbol{\pi},\boldsymbol{\mathfrak{H}})$ of
    $\boldsymbol{\mathcal{A}}$ such that
    $(\pi, \mathfrak{H})$ is unitarily equivalent to
    $\mathfrak{C}(\boldsymbol{\pi},\boldsymbol{\mathfrak{H}})$.
\end{definition}
Turning this definition into a question we arrive at the problem of
whether a given $^*$-representation of $\mathcal{A}$ can be
\emph{deformed} into a $^*$-representation of
$\boldsymbol{\mathcal{A}}$. In general this problem seems to be quite
hard to solve. In particular, some easy examples like the
Bargmann-Fock representation show that it is \emph{not} enough to
consider deformations with 
$\boldsymbol{\mathfrak{H}} = \mathfrak{H}[[\lambda]]$, see
\cite{waldmann:2000:pre}.
\begin{definition}
    \label{definition:posdef}
    Let $\omega_0: \mathcal{A} \to \ring{C}$ be a positive
    $\ring{C}$-linear functional. Then the
    $\ring{C}[[\lambda]]$-linear functional
    $\boldsymbol{\omega} 
    = \sum_{r=0}^\infty \lambda^r \omega_r : 
    \boldsymbol{\mathcal{A}} \to \ring{C}[[\lambda]]$
    is called a \emph{deformation} of $\omega_0$ if
    $\boldsymbol{\omega}$ is positive with respect to $\star$. If any
    positive linear functional of $\mathcal{A}$ can be deformed then
    $\boldsymbol{\mathcal{A}}$ is called a 
    \emph{positive deformation}. 
\end{definition}
The physical interpretation of this definition is that 
\emph{every classical state is the classical limit of a quantum
  state}. Obviously, from the physical intuition this should be a
feature of any reasonable quantization scheme. In fact, Hermitian star
products on symplectic manifolds are always positive deformations
\cite{bursztyn.waldmann:2000a}. From \cite{waldmann:2000:pre} we have:
\begin{lemma}
    Let 
    $\boldsymbol{\omega}: 
    \boldsymbol{\mathcal{A}} \to \ring{C}[[\lambda]]$ 
    be a deformation of a positive $\ring{C}$-linear functional
    $\omega_0: \mathcal{A} \to \ring{C}$. Then the GNS representation
    $(\boldsymbol{\pi_\omega}, \boldsymbol{\mathfrak{H}_\omega})$ is a
    deformation of the GNS representation 
    $(\pi_{\omega_0}, \mathfrak{H}_{\omega_0})$.
\end{lemma}
From this fact one easily obtains that a Hermitian deformation
$\boldsymbol{\mathcal{A}}$ of a unital $^*$-algebra $\mathcal{A}$ is
positive if and only if any GNS representation of $\mathcal{A}$ can be
deformed. In particular, if 
$\mathfrak{C}\boldsymbol{\Omega} = \psi_{\Unit}$ for 
some $\boldsymbol{\Omega} \in \boldsymbol{\mathfrak{H}}$ with
$\mathfrak{C}(\boldsymbol{\pi}, \boldsymbol{\mathfrak{H}}) =
(\pi_{\omega_0}, \mathfrak{H}_{\omega_0})$
then 
$\boldsymbol{\omega}(A) = 
\SP{\boldsymbol{\Omega}, \boldsymbol{\pi}(A)\boldsymbol{\Omega}}$ 
is a deformation of $\omega_0$. This again emphasizes the importance
of GNS representations.  
\section{Rieffel Induction}
\label{sec:rieffel}
Given two $^*$-algebras $\mathcal{A}$, $\mathcal{B}$ over $\ring{C}$
we want to learn something about their relation by comparing
$\Rep(\mathcal{A})$ and $\Rep(\mathcal{B})$ in a functorial way. This
wish can be made precise by considering a $\mathcal{B}$-$\mathcal{A}$
bimodule $\mathcal{E}$. As usual $\mathcal{E}$ shall also possess an
underlying $\ring{C}$-module structure. If $\mathcal{A}$ and
$\mathcal{B}$ are unital we always assume that the units act as units
on $\mathcal{E}$. In the non-unital case we require that $\mathcal{E}$
is strongly non-degenerate for both actions.

Given such a bimodule $\mathcal{E}$ and given a $^*$-representation
$(\pi, \mathfrak{H})$ of $\mathcal{A}$,
\begin{equation}
    \label{eq:ktilde}
    \widetilde{\mathfrak{K}}
    = \mathcal{E} \otimes_{\mathcal{A}} \mathfrak{H} 
\end{equation}
clearly becomes a $\mathcal{B}$-left module. The main idea of
\emph{Rieffel induction}, as known from $C^*$-algebra theory
\cite{landsman:1998,lance:1995}, is to impose additional structures on
$\mathcal{E}$ such that \eqref{eq:ktilde} eventually gives a
$^*$-representation of $\mathcal{B}$. 
We shall now sketch how the construction for $C^*$-algebras can be
extended to this purely algebraic framework of $^*$-algebras over
$\ring{C}$.
\begin{definition}
    \label{definition:ainnerprod}
    A \emph{$\mathcal{A}$-valued, $\mathcal{A}$-right linear Hermitian
      inner product}
    on a $\mathcal{A}$-right module $\mathcal{E}$ is a map
    $\SP{\cdot,\cdot}: \mathcal{E} \times \mathcal{E} \to \mathcal{A}$
    such that for all $x, y \in \mathcal{E}$ and $A \in \mathcal{A}$
    \begin{enumerate}
    \item $\SP{\cdot,\cdot}$ is $\ring C$-linear in the second
        argument,
    \item $\SP{x, y} = \SP{y, x}^*$,
    \item $\SP{x, y \cdot A} = \SP{x, y} A$, 
    \item $\SP{x, x} \in \mathcal{A}^+$.
    \end{enumerate}
    If $\mathcal{E}$ is even a $\mathcal{B}$-$\mathcal{A}$
    bimodule then $\SP{\cdot,\cdot}$ is \emph{compatible} with the
    $\mathcal{B}$-left action if in addition
    \begin{enumerate}
        \addtocounter{enumi}{4}
    \item $\SP{B \cdot x, y} = \SP{x, B^* \cdot y}$.
    \end{enumerate}
    Moreover, $\SP{\cdot,\cdot}$ satisfies the 
    \emph{positivity condition} \textbf{P} if
    \begin{equation}
        \label{eq:POSITIVE}
        \SP{x \otimes \phi, y \otimes \psi} :=
        \SP{\phi, \pi(\SP{x,y})\psi}
        \quad
        \textrm{for}
        \quad
        x \otimes \phi, y \otimes \psi \in \widetilde{\mathfrak{K}}
    \end{equation}
    extends to a positive semi-definite Hermitian inner product on
    $\widetilde{\mathfrak{K}}$ 
    for \emph{all} $^*$-repre\-sentations $(\pi, \mathfrak{H})$ of
    $\mathcal{A}$.
\end{definition}

Actually, the positivity condition \textbf{P} turns out to be rather
mild. It follows `almost' from $\SP{x,x} \in \mathcal{A}^+$. In case
of a $C^*$-algebra \eqref{eq:POSITIVE} is necessarily positive
semi-definite and one can give various sufficient conditions in the
general case, too, see \cite{bursztyn.waldmann:2001a} for a detailed
discussion.

Now let us assume we have a $\mathcal{B}$-$\mathcal{A}$ bimodule with
compatible $\mathcal{A}$-valued $\mathcal{A}$-right linear inner
product $\SP{\cdot,\cdot}$ such that \textbf{P} holds. Then we
consider the quotient
\begin{equation}
    \label{eq:Kdef}
    \mathfrak{K} = 
    \widetilde{\mathfrak{K}} 
    \; \big/ \; \widetilde{\mathfrak{K}}^\bot 
\end{equation}
of $\widetilde{\mathfrak{K}}$ divided by the vectors of length zero
and obtain a pre-Hilbert space $\mathfrak{K}$ over
$\ring{C}$. Moreover, we observe that the $\mathcal{B}$-left action on
$\widetilde{\mathfrak{K}}$ passes to the quotient $\mathfrak{K}$ and
yields a $^*$-representation of $\mathcal{B}$, thanks to the
compatibility condition. A final check shows that this construction is
functorial: indeed, let $(\pi_1, \mathfrak{H}_1)$ and 
$(\pi_2, \mathfrak{H}_2)$ be $^*$-representations of $\mathcal{A}$ and
let $U: \mathfrak{H}_1 \to \mathfrak{H}_2$ be an isometric
intertwiner. Then
\begin{equation}
    \label{eq:NewIntertwiner}
    \widetilde{V} (x \otimes \phi) = x \otimes U(\phi)
\end{equation}
passes to the quotient and yields an isometric intertwiner
$V: \mathfrak{K}_1 \to \mathfrak{K}_2$, Thus we end up with a functor
\begin{equation}
    \label{eq:rieffelind}
    \mathsf{R}_{\mathcal{E}}: 
    \Rep(\mathcal{A}) \longrightarrow \Rep(\mathcal{B})
\end{equation}
called \emph{Rieffel induction} functor. This functor turns out to be
compatible with algebraic constructions like direct sums
of $^*$-representations, tensor products of $^*$-algebras etc., see
\cite[Sect.~4]{bursztyn.waldmann:2001a}. In particular, Rieffel
induction can be viewed as a generalization of the GNS construction:
take $\mathcal{A}$ as $\mathcal{A}$-left module and as
$\ring{C}$-right module and set 
$\SP{A, B}_\omega := \omega(A^*B)$ for a positive linear functional
$\omega$ of $\mathcal{A}$. Then $\SP{\cdot,\cdot}_\omega$ satisfies
all needed requirements and $\textsf{R}_{\mathcal{A}}$ applied to the
canonical $^*$-representation of $\ring{C}$ on 
$\mathfrak{H} = \ring{C}$ gives the GNS representation 
$(\pi_\omega, \mathfrak{H}_\omega)$ of $\mathcal{A}$.
\section{Morita Equivalence of $^*$-Algebras}
\label{sec:morita}
While Rieffel induction is already an interesting and powerful tool
to understand the relation between $\Rep(\mathcal{A})$ and
$\Rep(\mathcal{B})$ it remains a `one way ticket'. We shall now
discuss whether and how one can find the `way back', i.e.\  from
$\Rep(\mathcal{B})$ to $\Rep(\mathcal{A})$, by imposing additional
structures on the bimodule $\mathcal{E}$. This will eventually lead to
a notion of Morita equivalence for $^*$-algebras over $\ring{C}$.

First we want a `symmetric situation' in $\mathcal{A}$ and
$\mathcal{B}$ whence we also require a $\mathcal{B}$-valued Hermitian
inner product
$\Theta_{\cdot,\cdot}: \mathcal{E} \times \mathcal{E} \to \mathcal{B}$
on $\mathcal{E}$. Since $\mathcal{E}$ is a $\mathcal{B}$-left module
we require it to be $\ring{C}$-linear as well as $\mathcal{B}$-left
linear in the \emph{first} argument. Beside this, the same
requirements as in Def.~\ref{definition:hermdef} are assumed to be
fulfilled. In particular we have now the compatibility 
\begin{equation}
    \label{eq:Acomp}
    \Theta_{x, y \cdot A} = \Theta_{x \cdot A^*, y}
\end{equation}
with the $\mathcal{A}$-right action on $\mathcal{E}$. Moreover, we say
$\SP{\cdot,\cdot}$ and $\Theta_{\cdot,\cdot}$ are \emph{compatible} if
\begin{equation}
    \label{eq:SPThetaComp}
    \Theta_{x,y} \cdot z = x \cdot \SP{y,z}
\end{equation}
for all $x, y \in \mathcal{E}$. Finally, the inner products are called
\emph{full} if
\begin{equation}
    \label{eq:Fullness}
    \mathcal{A} 
    = \ring{C}\textrm{-span } \{\SP{x,y} \; | \; x,y \in \mathcal{E}\}
    \quad
    \textrm{and}
    \quad
    \mathcal{B} 
    = \ring{C}\textrm{-span } \{\Theta_{x,y} \; | \; x,y \in \mathcal{E}\}.    
\end{equation}
Having all these requirements fulfilled we call $\mathcal{A}$ and
$\mathcal{B}$ strongly Morita equivalent as $^*$-algebras over
$\ring{C}$:
\begin{definition}
    \label{definition:morita}
    A $\mathcal{B}$-$\mathcal{A}$ bimodule $\mathcal{E}$ with
    compatible and full Hermitian inner products $\SP{\cdot,\cdot}$
    and $\Theta_{\cdot,\cdot}$ satisfying the above positivity
    requirements is called a \emph{Morita equivalence bimodule}. If
    there exists a Morita equivalence bimodule then $\mathcal{A}$ and
    $\mathcal{B}$ are called 
    \emph{strongly Morita equivalent as $^*$-algebras over}
    $\ring{C}$.
\end{definition}

Before we discuss some consequences of strong Morita equivalence we
shall point out that we indeed have a symmetric relation. Denote by
$\cc{\mathcal{E}}$ the complex-conjugate module, i.e.\ 
$\cc{\mathcal{E}} = \mathcal{E}$ as additive group but the
$\ring{C}$-module structure is given by $a \cc{x} = \cc{\cc{a} x}$,
where the identity map $\mathcal{E} \leftrightarrow \cc{\mathcal{E}}$
is denoted by $x \mapsto \cc{x}$. Clearly $\cc{\mathcal{E}}$ becomes a
$\mathcal{A}$-$\mathcal{B}$ bimodule by
\begin{equation}
    \label{eq:ccBimod}
    A \cdot \cc{x} = \cc{x \cdot A^*}
    \quad
    \textrm{and}
    \quad
    \cc{x} \cdot B = \cc{B^* \cdot x}.
\end{equation}
Similarly one obtains inner products with values in $\mathcal{A}$ and
$\mathcal{B}$, respectively. In case $\mathcal{A}$ and $\mathcal{B}$
are unital one gets indeed an equivalence of $\Rep(\mathcal{A})$ and
$\Rep(\mathcal{B})$, see \cite[Thm.~5.10]{bursztyn.waldmann:2001a}.
\begin{theorem}
    \label{theorem:morita}
    Let $\mathcal{A}$, $\mathcal{B}$ be unital $^*$-algebras over
    $\ring{C}$ and $\mathcal{E}$ a Morita equivalence bimodule.
    \begin{enumerate}
    \item $\cc{\mathcal{E}}$ is a Morita equivalence bimodule as well
        and 
        \begin{equation}
            \label{eq:equivCat}
            \Rep(\mathcal{A})
            \stackrel{\mathsf{R}_{\mathcal{E}}}{\longrightarrow}
            \Rep(\mathcal{B})
            \stackrel{\mathsf{R}_{\cc{\mathcal{E}}}}{\longrightarrow}
            \Rep(\mathcal{A})          
        \end{equation}
        gives an equivalence of categories.
    \item $\mathcal{A}$ and $\mathcal{B}$ are Morita equivalent as
        unital rings. In particular, $\mathcal{E}$ is finitely
        generated and projective as $\mathcal{A}$-right module and as
        $\mathcal{B}$-left module. 
    \end{enumerate}
\end{theorem}
The second statement has several important consequences: all the
`Morita invariants' known from Morita equivalence of unital rings are
also `strongly Morita invariant', see \cite{lam:1999} for a detailed
discussion of the Morita invariants. On the other hand the
$^*$-involution and the positivity structures imply more invariants,
see \cite[Thm.~5.4]{bursztyn.waldmann:2000a:pre}.
\begin{proposition}
    \label{proposition:lattice}
    Let $\mathcal{A}$, $\mathcal{B}$ be strongly Morita equivalent
    $^*$-algebras over $\ring{C}$. Then they have isomorphic lattices
    of closed ideals and 
    $\mathcal{J}_{\mathrm{min}}(\mathcal{A}) = \{0\}$ 
    if and only if
    $\mathcal{J}_{\mathrm{min}}(\mathcal{B}) = \{0\}$.
\end{proposition}
The last statement shows that the equivalence of the categories
$\Rep(\mathcal{A})$ and $\Rep(\mathcal{B})$ is in general a strictly
weaker notion than strong Morita equivalence as $^*$-algebras over
$\ring{C}$ since e.g.\  $\Rep(\mathcal{A})$ and
$\Rep(\mathcal{A}\big/\mathcal{J}_{\mathrm{min}}(\mathcal{A}))$ are
equivalent categories.
\begin{remark}
    \label{remark:morita}
    \hfill
    \begin{enumerate}
    \item Forgetting the order structure of $\ring{R}$ and thus all
        positivity requirements but keeping the $^*$-involution one
        arrives at Ara's notion of \emph{Morita $^*$-equivalence} for
        rings with involution, see \cite{ara:1999}. 
    \item As usual in Morita theory, the non-unital case is more
        involved. Here one obtains still reasonable results if the
        $^*$-algebras have an approximate identity. 
    \item In case of $C^*$-algebras one has the following relation
        between the different notions of Morita equivalence: two
        $C^*$-algebras are strongly Morita equivalent as
        $C^*$-algebras in Rieffel's sense\cite{rieffel:1974b} if and
        only if their \emph{Pedersen ideals} are strongly Morita
        equivalent as $^*$-algebras over $\ring{C}$, see
        \cite[Thm.~3.7]{bursztyn.waldmann:2000a:pre}.
    \end{enumerate}
\end{remark}

Let us now discuss a few examples for the unital case. Since here
$\mathcal{E}$ has to be projective we consider a \emph{Hermitian}
projection $P_0 \in M_n(\mathcal{A})$, i.e.\  $P_0^2 = P_0 = P_0^*$
and $\mathcal{E} = P_0 \mathcal{A}^n \subseteq \mathcal{A}^n$ as
$\mathcal{A}$-right module. In principle, this is an extra
requirement. But for $C^*$-algebras as well as for $C^\infty(M)$ any
projection turns out to be equivalent to a Hermitian projection. Thus
we shall always assume this in the following. Then we have
the canonical $\mathcal{A}$-valued inner product
\begin{equation}
    \label{eq:EinnerProd}
    \SP{x,y} = \sum_{i=1}^n x_i^* y_i
\end{equation}
for $x, y \in \mathcal{E}$, which is the restriction of the canonical
$\mathcal{A}$-valued $\mathcal{A}$-right linear Hermitian inner
product on the free $\mathcal{A}$-module $\mathcal{A}^n$. Next,
consider $\mathcal{B} = \End_{\mathcal{A}} (\mathcal{E}) 
= P_0 M_n (\mathcal{A}) P_0$,
which clearly is a $^*$-algebra acting on $\mathcal{E}$ such that
$\mathcal{E}$  becomes a $\mathcal{B}$-$\mathcal{A}$
bimodule. Moreover, this left action is clearly compatible with
\eqref{eq:EinnerProd}. Finally, we have the canonical
$\mathcal{B}$-valued inner product 
$\Theta_{x,y} \in \mathcal{B}$ \emph{defined} by
\begin{equation}
    \label{eq:BinnerProd}
    \Theta_{x,y} \cdot z = x \cdot \SP{y,z}.
\end{equation}
Thus the compatibility conditions will all be satisfied
automatically. However, in general $\mathcal{E}$ will be not a Morita
equivalence bimodule since the fullness conditions and the positivity
conditions are not fulfilled in general. Nevertheless there is a
simple \emph{sufficient} condition, see 
\cite[Def.~4.1 and Thm.~4.2]{bursztyn.waldmann:2000b}:
\begin{definition}
    \label{definition:StronglyFull}
    A Hermitian projection $P_0 \in M_n(\mathcal{A})$ is called
    \emph{strongly full} if there exists an invertible
    $\tau \in \mathcal{A}$ such that $\tr P_0 = \tau^*\tau$.
\end{definition}
\begin{theorem}
    \label{theorem:StrongylFull}
    Let $P_0 \in M_n(\mathcal{A})$ be a strongly full Hermitian
    projection. Then $\mathcal{E} = P_0\mathcal{A}^n$ is a Morita
    equivalence bimodule for 
    $\End_{\mathcal{A}}(\mathcal{E}) = P_0 M_n (\mathcal{A}) P_0$ and
    $\mathcal{A}$. 
\end{theorem}

Let us now discuss a more concrete example bringing us back to
deformation quantization: If $E \to M$ is a complex vector
bundle of fiber dimension  $k > 0$ with Hermitian fiber metric $h_0$
then the classical theorem of Serre and Swan, see
e.g.~\cite{swan:1962}, says that $\mathcal{E} = \Gamma^\infty(E)$ is
isomorphic to some $P_0 C^\infty(M)^n$ as
$C^\infty(M)$-modules. Moreover, since any two 
Hermitian fiber metrics are isometric, we can assume $P_0$ to be
Hermitian. In this case  $\tr P_0$ equals the fiber dimension of $E$
and thus $P_0$ is strongly full whence $\mathcal{E}$ is a Morita
equivalence bimodule for 
$\End_{C^\infty(M)}(\Gamma^\infty(E)) = \Gamma^\infty(\End(E))$ and
$C^\infty(M)$, see also \cite[Prop.~6.7]{bursztyn.waldmann:2000a:pre}
for a alternative proof.
\begin{lemma}
    \label{lemma:EndEME}
    Let $E \to M$ be a non-zero Hermitian vector bundle over $M$. Then
    $\Gamma^\infty(\End(E))$ and $C^\infty(M)$ are strongly Morita
    equivalent as $^*$-algebras over $\mathbb{C}$ via the equivalence
    bimodule $\Gamma^\infty(E)$. In particular, writing
    $\Gamma^\infty(E)$ as $P_0 C^\infty(M)^n$ the Hermitian projection
    $P_0$ is strongly full.
\end{lemma}
\section{Deformation of Projective Modules}
\label{sec:module}
Since isomorphism classes of Hermitian vector bundles $E \to M$ are in
one-to-one correspondence to isomorphism classes of Morita equivalence
bimodules for $C^\infty(M)$ we shall now deform the module structure
of $\Gamma^\infty(E)$ as well as the fiber metric $h_0$ in order to
obtain equivalence bimodules for star products. Note that there are
other motivations for deforming the (Hermitian) vector bundle
structure like e.g.\  the index theorems of deformation quantization
\cite{nest.tsygan:1995a,fedosov:1996} or the global description of
non-commutative field theories arising from string theory
\cite{jurco.schupp.wess:2001a:pre,jurco.schupp.wess:2001b:pre,jurco.schupp.wess:2000,waldmann:2001a:pre}.

Again we consider the general situation where $\mathcal{A}$ is a
unital $^*$-algebra over $\ring{C}$ and 
$\boldsymbol{\mathcal{A}} = (\mathcal{A}[[\lambda]],\star)$ is a
Hermitian deformation of $\mathcal{A}$. Here and in the following we
shall assume that $\mathbb{Q} \subseteq \ring{R}$. Let 
$P_0 \in M_n (\mathcal{A})$ be a Hermitian projection and
$\mathcal{E} = P_0 \mathcal{A}^n$ as before. Since the matrices
$M_n(\boldsymbol{\mathcal{A}})$ inherit the deformed product $\star$
we can consider
\begin{equation}
    \label{eq:defproj}
    \boldsymbol{P} = 
    \frac{1}{2} + \left(P_0 - \frac{1}{2}\right) \star
    \frac{1}{\sqrt[\star]{1 + 4 (P_0 \star P_0 - P_0)}},
\end{equation}
which gives a well-defined formal power series 
$\boldsymbol{P} \in M_n(\boldsymbol{\mathcal{A}})$. It is an easy
computation to show that 
$\boldsymbol{P} \star \boldsymbol{P} 
= \boldsymbol{P} = \boldsymbol{P}^*$
is a Hermitian projection with respect to $\star$, see
\cite[Eq.~(6.1.4)]{fedosov:1996} or \cite{gerstenhaber.schack:1988}.
Thus one always can deform Hermitian projections since clearly
$\boldsymbol{P}$ coincides with $P_0$ in zeroth order.
Moreover,
\begin{equation}
    \label{eq:IsoProj}
    I: \mathcal{E}[[\lambda]] \ni x 
    \; \mapsto \;
    \boldsymbol{P} \star x \in 
    \boldsymbol{P} \star \boldsymbol{\mathcal{A}}^n 
\end{equation}
turns out to be an isomorphism of $\ring{C}[[\lambda]]$-modules. Since
the right hand side is obviously a projective $\star$-right module we
can pull-back this module structure via $I$ to obtain a finitely
generated projective $\star$-right module 
$(\mathcal{E}[[\lambda]], \bullet)$ which is a \emph{deformation} of
the original module $\mathcal{E}$. Moreover, the restriction of the
canonical $\boldsymbol{\mathcal{A}}$-valued
$\boldsymbol{\mathcal{A}}$-right linear Hermitian inner product on 
$\boldsymbol{\mathcal{A}}^n$ to 
$\boldsymbol{P} \star \boldsymbol{\mathcal{A}}^n$ induces via
\eqref{eq:IsoProj} a Hermitian inner product $\boldsymbol{h}$ on
$(\mathcal{E}[[\lambda]],\bullet)$, which is again a deformation of the
original $h_0$. Analogously to \eqref{eq:IsoProj} the map
\begin{equation}
    \label{eq:IsoProjEndos}
    P_0 M_n(\mathcal{A}[[\lambda]])P_0 \ni B 
    \; \mapsto \; 
    \boldsymbol{P} \star B \star \boldsymbol{P} \in 
    \boldsymbol{P} \star M_n (\boldsymbol{\mathcal{A}}) 
    \star \boldsymbol{P} 
    \cong \End\nolimits_{\boldsymbol{\mathcal{A}}} 
    (\boldsymbol{P} \star \boldsymbol{\mathcal{A}}^n)
\end{equation}
is an isomorphism of $\ring{C}[[\lambda]]$-modules. Hence the
pull-back of the operator product to the left hand side gives a
deformation $\star'$ of $\End_{\mathcal{A}} (\mathcal{E})$. 
Since \eqref{eq:IsoProjEndos} is compatible with the $^*$-involutions,
$\star'$ is Hermitian whence
$(\End_{\mathcal{A}}(\mathcal{E})[[\lambda]], \star') \cong
\End_{\boldsymbol{\mathcal{A}}} 
(\mathcal{E}[[\lambda]], \bullet, \boldsymbol{h})$ as $^*$-algebras.
Clearly, any choice of an isomorphism of the algebras
$(\End_{\mathcal{A}}(\mathcal{E})[[\lambda]],\star')$ and
$\End_{\boldsymbol{\mathcal{A}}}(\mathcal{E}[[\lambda]],
\bullet,\boldsymbol{h})$
gives automatically a left action of the first algebra on
$\mathcal{E}[[\lambda]]$ by $\bullet$-right linear endomorphisms. 
Using the particular isomorphism \eqref{eq:IsoProjEndos} this left
action is obviously a \emph{deformation} of the original left action
of $\End_{\mathcal{A}}(\mathcal{E})$ on $\mathcal{E}$.

Finally we remark that any other deformations 
$(\bullet', \boldsymbol{h}')$ of the original module structure and
inner product are necessarily equivalent to the ones above. This is a
direct consequence of $(\mathcal{E}[[\lambda]], \bullet)$ being
\emph{projective} over $\boldsymbol{\mathcal{A}}$ as well as 
Lemma~\ref{lemma:deformsquares}. Hence all choices lead to isomorphic
deformations $\star'$. The above choice is distinguished in so far
that it gives a deformation of the bimodule structure of
$\mathcal{E}$. It can be shown easily that any other $\star'$
admitting a deformation of the left action is actually
\emph{equivalent} to the above choice (and not only isomorphic). In
the following we shall always use such a $\star'$.
\begin{theorem}
    \label{theorem:defprojmod}
    Let $\mathcal{A}$ be a $^*$-algebra over $\ring{C}$ and 
    $P_0 = P_0^2 = P_0^* \in M_n(\mathcal{A})$ and
    $\mathcal{E} = P_0\mathcal{A}^n$ endowed with its
    $\mathcal{A}$-valued $\mathcal{A}$-right linear Hermitian inner
    product $h_0$ inherited from $\mathcal{A}^n$. Moreover, let
    $\boldsymbol{\mathcal{A}} = (\mathcal{A}[[\lambda]],\star)$ be a
    Hermitian deformation of $\mathcal{A}$. Then there exists a
    deformation $\bullet$ of the right module structure of
    $\mathcal{E}$ with respect to $\star$ as well as a deformation
    $\boldsymbol{h}$ of $h_0$. Any two such deformations are
    (isometrically) equivalent. Moreover, $\bullet$ yields a finitely
    generated projective $\boldsymbol{\mathcal{A}}$-module. Finally,
    $\bullet$ induces a Hermitian deformation $\star'$ of 
    $\mathcal{B} = \End_{\mathcal{A}}(\mathcal{E})$ such that
    $(\mathcal{B}[[\lambda]], \star') \cong
    \End_{\boldsymbol{\mathcal{A}}} 
    (\mathcal{E}[[\lambda]], \bullet)$ with a deformed bimodule
    structure, and $\star'$ is uniquely determined up to equivalence.
\end{theorem}
\begin{remark}
    \label{remark:defproj}
    \hfill
    \begin{enumerate}
    \item Since for any Hermitian vector bundle $E \to M$ the sections
        $\mathcal{E} = \Gamma^\infty(E)$ are exactly of this form we
        obtain the result that \emph{Hermitian vector bundles can
          always be deformed in a unique way up to equivalence}. 
        See also \cite{lecomte.roger:1988} for the case of a trivial
        vector bundle.
    \item Since clearly any projection 
        $\boldsymbol{P} \in M_n(\boldsymbol{\mathcal{A}})$ has a
        projection $P_0 \in M_n (\mathcal{A})$ as classical limit it
        follows that all finitely generated projective modules over
        $\boldsymbol{\mathcal{A}}$ are deformations of finitely
        generated projective modules over $\mathcal{A}$. Thus the
        \emph{$K$-theory does not change under deformation}, i.e.\
        $K_0(\mathcal{A}) \cong K_0 (\boldsymbol{\mathcal{A}})$, see
        e.g.~\cite{rosenberg:1996:pre}.
    \end{enumerate}
\end{remark}

The use of strongly full projections is illustrated by the following
lemma, which follows directly form Lemma~\ref{lemma:deformsquares}.
\begin{lemma}
    \label{lemma:deffull}
    Let $P_0 \in M_n(\mathcal{A})$ be a strongly full projection then
    any deformation 
    $\boldsymbol{P} \in M_n (\boldsymbol{\mathcal{A}})$ is again
    strongly full.
\end{lemma}
Thus if $\mathcal{E} = P_0 \mathcal{A}^n$ is a Morita equivalence
bimodule for $\mathcal{A}$ and 
$\mathcal{B} = \End_{\mathcal{A}} (\mathcal{E})$ with a strongly full
$P_0$ then $\boldsymbol{\mathcal{A}}$ and
$\boldsymbol{\mathcal{B}} = (\mathcal{B}[[\lambda]],\star')$, with
$\star'$ as above, are again strongly Morita equivalent as
$^*$-algebras over $\ring{C}[[\lambda]]$ via
$(\mathcal{E}[[\lambda]],\bullet,\boldsymbol{h})$ by
Theorem~\ref{theorem:defprojmod}.
\section{Local Description of Deformed Vector Bundles}
\label{sec:vector}
Before we shall investigate the Morita equivalence of star products we
shall give a local description of the deformed Hermitian vector
bundles in case where $\mathcal{A} = C^\infty(M)$ and
$\mathcal{E} = \Gamma^\infty(E)$, see also
\cite{jurco.schupp.wess:2001a:pre} for the case of a line bundle. Let
$\star$ be a Hermitian star product for $M$ and $\bullet$,
$\boldsymbol{h}$ some fixed choices of deformations for the Hermitian
vector bundle. In particular, $\bullet$ and $\boldsymbol{h}$ can be
chosen to be \emph{local} thanks to the explicit map
\eqref{eq:IsoProj}.

Now we choose a good open and locally finite cover
$\{\mathcal{O}_\alpha\}$ of $M$ and denote by  
$e_\alpha = (e_{\alpha,1},\ldots,e_{\alpha,k})$ a local
frame on $\mathcal{O}_\alpha$ of the undeformed vector bundle, i.e.\
the $e_{\alpha,i}$ give local base sections and any local section $s$
can be written as $s = \sum_{i=1}^k e_{\alpha,i} s_\alpha^i$
with unique coefficient functions 
$s_\alpha^i \in C^\infty(\mathcal{O}_\alpha)$. For abbreviation we
simply write $s_\alpha = (s_\alpha^1, \ldots, s_\alpha^k)$ and 
$s = e_\alpha s_\alpha$. A simple induction shows the quantum analog
of this:
\begin{lemma}
    \label{lemma:qframe}
    $\qe_{\alpha,i} 
    = e_{\alpha,i} + \lambda e_{\alpha,i}^{(1)} + \cdots 
    \in \Gamma^\infty(E)[[\lambda]]$
    be deformations of the classical base sections
    $e_{\alpha,i}$. Then any local section 
    $\qs \in \Gamma^\infty(E|_{\mathcal{O}_\alpha})[[\lambda]]$ can
    uniquely be written as
    \begin{equation}
        \label{eq:qscoeff}
        \qs = \sum_{i=1}^k \qe_{\alpha,i} \bullet \qs_\alpha^i
    \end{equation}
    with 
    $\qs_\alpha^i \in C^\infty(\mathcal{O}_\alpha)[[\lambda]]$.
    Moreover, the $\mathbb{C}[[\lambda]]$-linear map 
    $\qpsi_\alpha: \Gamma^\infty(E|_{\mathcal{O}_\alpha})[[\lambda]]
    \to C^\infty (\mathcal{O}_\alpha)^k[[\lambda]]$ defined by
    \begin{equation}
        \label{eq:qtriv}
        \qpsi_\alpha(\qs) = \qs_\alpha 
        = (\qs_\alpha^1, \ldots, \qs_\alpha^k)
    \end{equation}
    is an isomorphism of $\star$-right modules, i.e. 
    $\qpsi_\alpha(\qs \bullet f) = \qpsi_\alpha(\qs) \star f$.
\end{lemma}
Clearly the $\qpsi_\alpha$ can be viewed as deformations of the
classical trivialization isomorphisms $\psi_\alpha$ determined by the
classical frame $e_\alpha$. We shall call 
$\qe_\alpha = (\qe_{\alpha,1}, \ldots, \qe_{\alpha,k})$ a local frame
for the deformed vector bundle and write 
$\qs = \qe_\alpha \bullet \qs_\alpha$ for short. In the case where we
have also deformed the Hermitian fiber metric we can apply again
Lemma~\ref{lemma:deformsquares} to obtain orthonormal frames:
\begin{lemma}
    \label{lemma:orthoframe}
    Let $\widetilde{\qe}_\alpha$ be a local frame such that the zeroth
    orders constitute an orthonormal frame with respect to $h_0$. Then
    there exists a matrix 
    $\qV = \Unit + \sum_{r=1}^\infty \lambda^r \qV_r \in
    M_k (C^\infty(\mathcal{O}_\alpha))[[\lambda]]$ such that 
    $\qe_\alpha = \widetilde{\qe}_\alpha \bullet \qV$ is an
    orthonormal frame with respect to $\qh$, i.e.
    \begin{equation}
        \label{eq:orthoframe}
        \qh(\qe_{\alpha,i}, \qe_{\alpha,j}) = \delta_{ij}.
    \end{equation}
\end{lemma}
As in the undeformed case the choice of local frames gives local
\emph{transition matrices} which satisfy the cocycle identity:
\begin{lemma}
    \label{lemma:cocycle}
    Let $\{\qe_\alpha\}$ be local frames then there exist unique
    transition matrices 
    $\qphi_{\alpha\beta} \in 
    M_k(C^\infty(\mathcal{O}_\alpha\cap\mathcal{O}_\beta))[[\lambda]]$
    on every non-empty overlap such that
    \begin{equation}
        \label{eq:qtrans}
        \qe_\beta = \qe_\alpha \bullet \qpsi_{\alpha\beta}
        \quad
        \textrm{and}
        \quad
        \qs_\alpha = \qpsi_{\alpha\beta} \star \qs_\beta
        = \qpsi_\alpha \circ \qpsi^{-1}_\beta (\qs_\beta).
    \end{equation}
    Moreover, the cocycle condition 
    \begin{equation}
        \label{eq:cocycle}
        \qphi_{\alpha\beta} \star \qphi_{\beta\alpha} = \Unit
        \quad
        \textrm{and}
        \quad
        \qphi_{\alpha\beta} \star 
        \qphi_{\beta\gamma} \star
        \qphi_{\gamma\alpha} = \Unit
    \end{equation}
    is fulfilled on the corresponding non-empty intersections. If
    the $\{\qe_\alpha\}$ are orthonormal we have unitary transition
    matrices $\qphi_{\alpha\beta}^* = \qphi_{\beta\alpha}$.    
\end{lemma}
Again the $\qphi_{\alpha\beta}$ are deformations of the transition
matrices $\phi_{\alpha\beta}$ corresponding to the classical limit of
the frames $\{\qe_\alpha\}$.

Let us now discuss the $\bullet$-right linear endomorphisms of the
deformed vector bundle. If $\qA (\qs \bullet f) = \qA(\qs) \bullet f$
holds globally then locally we have a matrix 
$\qA_\alpha \in M_k (C^\infty(\mathcal{O}_\alpha))[[\lambda]]$ such
that
\begin{equation}
    \label{eq:localendo}
    (\qA(\qs))_\alpha = \qA_\alpha \star \qs_\alpha
    \quad
    \textrm{and}
    \quad
    \qA(\qe_\alpha) = \qe_\alpha \bullet \qA_\alpha.
\end{equation}
On overlaps $\mathcal{O}_\alpha \cap \mathcal{O}_\beta$ we clearly
have
\begin{equation}
    \label{eq:endooverlap}
    \qA_\beta = 
    \qphi_{\beta\alpha} \star \qA_\alpha \star \qphi_{\alpha\beta}.
\end{equation}
Conversely, if a collection of matrices $\{\qA_\alpha\}$ satisfies
\eqref{eq:endooverlap} then it defines a global $\bullet$-right linear
endomorphism $\qA$ by \eqref{eq:localendo}. The composition of $\qA$
and $\qB$ is locally given by the deformed matrix multiplication
\begin{equation}
    \label{eq:ABlocally}
    (\qA \circ \qB)_\alpha = \qA_\alpha \star \qB_\alpha,
\end{equation}
and if the frames $\{\qe_\alpha\}$ are orthonormal we have
$(\qA^*)_\alpha = (\qA_\alpha)^*$. 
%
% The following pagebreak makes the whole thing look much
% nicer. Perhaps it has to be changed when putting in the final form.
%

%\pagebreak

%
%
\section{From local to global deformations}
\label{sec:ltog}
In this section we present some new material, namely a detailed
comparison between the approach to deformed vector bundles by deformed
local transition matrices as used in
\cite{jurco.schupp.wess:2001a:pre,jurco.schupp.wess:2001b:pre} with
our approach of globally deforming a projection. As we already know
from Theorem~\ref{theorem:defprojmod} both approaches are necessarily
equivalent. The local approach is in some sense more familiar and more
suited for explicit constructions while in the global approach
questions of uniqueness up to equivalence can be handled more easily.

The starting point is an open cover $\{\mathcal{U}_\alpha\}$ of $M$
such that $E|_{\mathcal{U}_\alpha}$ is trivial. It will be necessary
to assume that $\{\mathcal{U}_\alpha\}$ is actually a \emph{finite}
cover whence we shall use a different symbol $\mathcal{U}_\alpha$
here, see e.g. \cite[Prop.~4.1, Chap.~III]{wells:1980} for the 
existence of such covers. We assume that we have found deformations
$\qphi_{\alpha\beta}$ of the classical transition matrices
$\phi_{\alpha\beta}$ corresponding to some local frames
$\{e_\alpha\}$, i.e.\  we have 
$\qphi_{\alpha\beta} = \phi_{\alpha\beta} + \sum_{r=1}^\infty
\lambda^r \qphi_{\alpha\beta}^{(r)}$ satisfying \eqref{eq:cocycle}. We
shall mainly be interested in the Hermitian case whence we assume in
addition $\qphi_{\alpha\beta} ^* = \qphi_{\beta\alpha}$. Then the aim
is to build a deformation $\bullet$ of $\Gamma^\infty(E)$ and $\star'$
of $\Gamma^\infty(\End(E))$ out of this data and to show that these
deformations are of the form as obtained in Section~\ref{sec:module}.

We choose a partition of unity $\{\chi_\alpha\}$ subordinate to
$\{\mathcal{U}_\alpha\}$. For $s \in \Gamma^\infty(E)[[\lambda]]$ we
define
\begin{equation}
    \label{eq:qslocal}
    \qs_\alpha := \sum_{\beta} \qphi_{\alpha\beta} \star \chi_\beta
    \star s_\beta,
\end{equation}
where classically $s = e_\beta s_\beta$ on
$\mathcal{U}_\beta$. Clearly 
$\qs_\alpha \in C^\infty(\mathcal{U}_\alpha)^k[[\lambda]]$ and
$\qphi_{\alpha\beta} \star \qs_\beta = \qs_\alpha$. Finally note that
$\qs_\alpha = s_\alpha + \sum_{r=1}^\infty \lambda^r
\qs_\alpha^{(r)}$.
Then the following lemma is a simple induction and can be seen as the
`converse' statement to Lemma~\ref{lemma:qframe}. 
\begin{lemma}
    \label{lemma:loctoglobsec}
    The map $s \mapsto \{\qs_\alpha\}$ from
    $\Gamma^\infty(E)[[\lambda]]$ into the set of all collections
    $\{\qs_\alpha\}$ with $\qs_\alpha \in
    C^\infty(\mathcal{U}_\alpha)^k[[\lambda]]$ satisfying
    $\qphi_{\alpha\beta} \star \qs_\beta = \qs_\alpha$ is a
    $\mathbb{C}[[\lambda]]$-linear bijection.
\end{lemma}
This allows to define a global $\star$-right module structure on the
sections by setting
\begin{equation}
    \label{eq:defbullet}
    (s \bullet f)_\alpha := \qs_\alpha \star f,
\end{equation}
which is clearly a deformation of the classical module structure. Thus
the existence of a deformation quantization of $E$ follows from
the existence of deformed transition matrices. A Hermitian inner
product is defined by
\begin{equation}
    \label{eq:defherm}
    \qh(s, t)|_{\mathcal{U}_\alpha} = \SP{\qs_\alpha, \qt_\alpha} =
    \sum_{l=1}^k \cc{\qs_\alpha^l} \star \qt_\alpha^l,
\end{equation}
which is globally defined since 
$\qphi_{\alpha\beta}^* = \qphi_{\beta\alpha}$ and clearly a
deformation of a classical Hermitian fiber metric $h_0$.

Now let us consider a `local' way to deal with the $\bullet$-right
linear endomorphisms. First we choose a \emph{quadratic} partition of
unity $\{\chi_\alpha\}$ subordinate to $\{\mathcal{U}_\alpha\}$ and
define for $A \in \Gamma^\infty (\End(E))[[\lambda]]$ the matrix
$\widetilde{T}_\alpha (A) 
\in M_k(C^\infty(\mathcal{U}_\alpha))[[\lambda]]$
by
\begin{equation}
    \label{eq:deftildeT}
    \widetilde{T}_\alpha(A) 
    = \sum_\gamma \qphi_{\alpha\gamma} \star \cc{\chi}_\gamma
    \star A_\gamma \star \chi_\gamma \star \qphi_{\gamma\alpha},
\end{equation}
where $A_\gamma$ is the local matrix of $A$ corresponding to the
undeformed structures. Then clearly
\begin{equation}
    \label{eq:widetildeT}
    \qphi_{\alpha\beta} \star \widetilde{T}_\beta (A) \star
    \qphi_{\beta\alpha} = \widetilde{T}_\alpha(A)
    \quad
    \textrm{and}
    \quad
    \widetilde{T}_\alpha(A^*) = \widetilde{T}_\alpha(A)^*,
\end{equation}
whence 
$(\widetilde{T}(A) s)_\alpha 
= \widetilde{T}_\alpha(A) \star \qs_\alpha$
defines a global $\bullet$-right linear endomorphism
$\widetilde{T}(A)$. In zeroth order we have 
$\widetilde{T}_\alpha(A) = A_\alpha$ whence a simple induction shows
that all $\bullet$-right linear endomorphisms are obtained this
way, i.e.\  $\widetilde{T}$ is a bijection. Thus we obtain a Hermitian
deformation $\widetilde{\star}$ of $\Gamma^\infty(\End(E))[[\lambda]]$
such that 
$\widetilde{T}_\alpha (A \,\widetilde{\star}\, B) =
\widetilde{T}_\alpha (A) \star \widetilde{T}_\alpha (B)$.
However, $\widetilde{T}_\alpha (\Unit) \ne \Unit$ in general. Thus the
identity endomorphism is no longer the unit for $\widetilde{\star}$
but some other
$E = \widetilde{T}^{-1}(\Unit) 
= \Unit + \sum_{r=1}^\infty \lambda^r E_r$.
Nevertheless this can be `repaired' in the usual way by passing to an
equivalent deformation $\star'$. By Lemma~\ref{lemma:deformsquares} we
know that $\Unit = (U^{-1})^* \,\widetilde{\star}\, U^{-1}$ with some
$\widetilde{\star}$-invertible $U$. The equivalence transformation
$S(A) := 
(U^{-1})^* \,\widetilde{\star}\, A \,\widetilde{\star}\, U^{-1}$
gives a new deformation 
$S(A \,\widetilde{\star}\, B) = S(A) \star' S(B)$,
which is still Hermitian and satisfies 
$\Unit \star' A = A = A \star'\Unit$. Moreover, we obtain a new
isomorphism $T(A) = \widetilde{T}(S^{-1}(A))$ locally given by
\begin{equation}
    \label{eq:localT}
    T_\alpha(A) 
    = \widetilde{T}_\alpha (U^*) \star
    \widetilde{T}_\alpha(A) \star \widetilde{T}_\alpha (U)
    = \sum_\gamma  
    \boldsymbol{U}^*_\alpha \star\qphi_{\alpha\gamma} \star
    \cc{\chi}_\gamma \star A_\gamma \star \chi_\gamma \star
    \qphi_{\gamma\alpha} \star \boldsymbol{U}_\alpha, 
\end{equation}
where $\boldsymbol{U}_\alpha := \widetilde{T}_\alpha (U)$. Now we
indeed have 
\begin{equation}
    \label{eq:niceT}
    T_\alpha (A \star' B) = T_\alpha (A) \star T_\alpha (B),
    \quad
    T_\alpha (A^*) = T_\alpha (A)^*,
    \quad
    \textrm{and}
    \quad
    T_\alpha (\Unit) = \Unit.
\end{equation}
Again $T$ as well as $\star'$ are not uniquely determined and there
may be situations where a more concrete description of the isomorphism
$T$ is available.

Let us now prove `intrinsically' that such a deformation $\bullet$
yields a finitely generated projective module. This can be viewed as a
sort of \emph{quantum Serre-Swan theorem}. We follow the ideas of
\cite{jurco.schupp.wess:2001b:pre} where the case of a line bundle was
considered. It is crucial to have a finite cover
$\{\mathcal{U}_\alpha\}$.

Let $\chi_\alpha$ be as before, only viewed as element in
$\Gamma^\infty(\End(E))$. 
Using Lemma~\ref{lemma:deformsquares} we find a $\star'$-invertible
$V \in \Gamma^\infty(\End(E))[[\lambda]]$ such that
\begin{equation}
    \label{eq:chichiVV}
    \sum_{\alpha = 1}^m \cc{\chi}_\alpha \star' \chi_\alpha 
    = V^* \star' V,
\end{equation}
Hence the global endomorphisms 
$\boldsymbol{\chi}_\alpha = \chi_\alpha \star' V^{-1}$ satisfy
\begin{equation}
    \label{eq:boldchi}
    \sum_{\alpha = 1}^m 
    \boldsymbol{\chi}_\alpha^* \star' \boldsymbol{\chi}_\alpha 
    = \Unit 
    \quad
    \textrm{and}
    \quad
    \supp{\boldsymbol{\chi}_\alpha} \subseteq \mathcal{U}_\alpha.
\end{equation}
We define a $\mathbb{C}[[\lambda]]$-linear map 
$\epsilon:
\Gamma^\infty(E)[[\lambda]] \to C^\infty(M)^n[[\lambda]]$,
where $n = mk$, by
\begin{equation}
    \label{eq:epsilondef}
    \epsilon(s) 
    = \left(T_1(\boldsymbol{\chi}_1) \star \qs_1, 
        \ldots, T_m(\boldsymbol{\chi}_m) \star \qs_m 
    \right),
\end{equation}
which indeed gives \emph{global} functions since $T$ and $\star$ are
local and 
$\supp\boldsymbol{\chi}_\alpha \subseteq \mathcal{U}_\alpha$. 
Clearly we have $\epsilon(s \bullet f) = \epsilon(s) \star f$,
where $C^\infty(M)^n[[\lambda]]$ is endowed with its usual \emph{free}
$\star$-right module structure. Next we define 
$\pi: C^\infty(M)^n[[\lambda]] \to \Gamma^\infty(E)[[\lambda]]$ by
\begin{equation}
    \label{eq:pidef}
    (\pi(t))_\alpha := \sum_{\gamma = 1}^m \qphi_{\alpha\gamma} \star
    T_\gamma (\boldsymbol{\chi}_\gamma^*) \star t_\gamma
    \big|_{\mathcal{U}_\alpha},
\end{equation}
where $t = (t_1, \ldots, t_m)$ with global 
$t_\gamma \in C^\infty(M)^k[[\lambda]]$. Then 
$\qphi_{\alpha\beta} \star \pi(t)_\beta = \pi(t)_\alpha$ whence
$\pi(t)$ is indeed a section of $E$. Again $\pi$ is a $\star$-right
module morphism and a straightforward computation yields $\pi \circ
\epsilon = \id$, whence $\epsilon \circ \pi$ is a projection. Since it
is even $\star$-right linear it is of the form
\begin{equation}
    \label{eq:epspiP}
    \epsilon(\pi(t)) = \boldsymbol{P} \star t
    \quad
    \textrm{with}
    \quad
    \boldsymbol{P} \star \boldsymbol{P} = \boldsymbol{P} 
    \in M_n (C^\infty(M))[[\lambda]].
\end{equation}
But this shows that $(\Gamma^\infty(E)[[\lambda]], \bullet)$ is indeed
a finitely generated projective module as $\epsilon$ provides a module
isomorphism into the image of $\boldsymbol{P}$. Thus we have an
alternative proof of this fact and can now easily use the uniqueness
statements of Theorem~\ref{theorem:defprojmod} also in this
(completely equivalent) approach.
\begin{example}[Flat bundles]
    \label{example:flatbundle}
    Let $E \to M$ be a \emph{flat} Hermitian vector bundle and let
    $\{e_\alpha\}$ be covariantly constant local frames with
    \emph{constant} transition functions $\phi_{\alpha\beta}$. Then
    $\qphi_{\alpha\beta} = \phi_{\alpha\beta}$ satisfy the cocycle
    condition \eqref{eq:cocycle}. Hence $\qs_\alpha = s_\alpha$ and 
    $s \bullet f = e_\alpha (s_\alpha \star f)$ if classically 
    $s = e_\alpha s_\alpha$. Moreover, in this case 
    $T_\alpha(A) = A_\alpha$ will do the job, i.e.\  
    \begin{equation}
        \label{eq:Tflat}
        (T(A)s)_\alpha = T_\alpha(A) \star s_\alpha = A_\alpha \star
        s_\alpha 
    \end{equation}
    gives a simple description of the $\bullet$-right linear
    endomorphisms and we have
    \begin{equation}
        \label{eq:starprimeflat}
        (A \star' B)_\alpha = A_\alpha \star A_\beta.
    \end{equation}
    Surprisingly, in this case the center $C^\infty(M)$ of the
    undeformed algebra $\Gamma^\infty(\End(E))$ is still a
    \emph{sub-algebra} with respect to $\star'$. Moreover, from
    \eqref{eq:starprimeflat} it immediately follows that $\star' =
    \star$ when restricted to 
    $C^\infty(M)[[\lambda]] \subseteq
    \Gamma^\infty(\End(E))[[\lambda]]$. 
    Thus a flat vector bundle becomes automatically a
    \emph{$\star$-bimodule}, see also Cor.~\ref{corollary:moritaII}.
\end{example}
It seems to be an interesting question whether in the case of a
general vector bundle the center is deformed into a sub-algebra with
respect to $\star'$. In this case $\star'$ would induce an
associative deformation of the center by restricting $\star'$.
\begin{question}
    Let $\mathcal{A}$ be an \emph{Abelian} algebra with associative
    deformation $\star$ and let $\mathcal{E}$ be a finitely generated
    projective module over $\mathcal{A}$ with deformation $\bullet$
    and let $\star'$ be the corresponding (unique up to equivalence)
    deformation of $\End_{\mathcal{A}}(\mathcal{E})$. 
    Then $\mathcal{A} \subseteq \End_{\mathcal{A}} (\mathcal{E})$
    belongs to the center. Is this deformed into a \emph{sub-algebra}
    with respect to $\star'$ or, if not, what is the algebra generated
    by $\mathcal{A}$ inside
    $\End_{\mathcal{A}}(\mathcal{E})[[\lambda]]$ with respect to
    $\star'$? Moreover, does the assumption of \emph{fullness} of
    $\mathcal{E}$ simplify the problem? 
\end{question}
\section{Morita Equivalence of Star Products}
\label{sec:mestar}
In this final section we shall sketch the classification of star
products up to Morita equivalence on symplectic manifolds.

First we recall the construction of the 
\emph{relative characteristic class} of two star products $\star$,
$\star'$ on $(M, \omega)$ as it can be found in
\cite{gutt.rawnsley:1999}. The \emph{characteristic class}
$c(\star)$ of a symplectic star product can e.g.\  be defined by means
of the Weyl curvature in the Fedosov construction, see
e.g.~\cite{nest.tsygan:1995a,bertelson.cahen.gutt:1997}, and yields a
class
\begin{equation}
    \label{eq:charclass}
    c(\star) \in \frac{1}{\im\lambda} [\omega] 
    + \HdR^2(M, \mathbb{C})[[\lambda]].
\end{equation}
Then the relative class of two star products is defined by
\begin{equation}
    \label{eq:relclass}
    t(\star',\star) = c(\star') - c(\star) \in 
    \HdR^2(M, \mathbb{C})[[\lambda]],
\end{equation}
and $\star'$ is equivalent to $\star$ if and only if their relative
class vanishes. There is a simple \v{C}ech cohomological construction
of $t(\star',\star)$. We again use the good cover as before. Then on
$\mathcal{O}_\alpha$ the star products $\star'$ and $\star$ are
equivalent via some 
$T_\alpha = \id + \sum_{r=1}^\infty \lambda^r T_\alpha^{(r)}$, i.e.\
$T_\alpha(f \star' g) = T_\alpha f \star T_\alpha g$, since
$\mathcal{O}_\alpha$ is contractible. On overlaps it follows that
$T_\alpha^{-1}T_\beta$ is a $\star$-automorphism starting with the
identity. Hence it is an inner automorphism
\begin{equation}
    \label{eq:localautos}
    T_\alpha^{-1} T_\beta (f) = \eu^{[t_{\alpha\beta}, \cdot]} (f)
    = \Exp(t_{\alpha\beta}) \star f \star \Exp(-t_{\alpha\beta})
\end{equation}
with some function 
$t_{\alpha\beta} 
\in C^\infty(\mathcal{O}_\alpha\cap\mathcal{O}_\beta)[[\lambda]]$. 
Here $\Exp(\cdot)$ denotes the $\star$-exponential, see
e.g.~\cite{bordemann.roemer.waldmann:1998,bursztyn.waldmann:2001:pre}.
Moreover, the $T_\alpha^{-1} T_\beta$ satisfy the `cocycle condition'
\begin{equation}
    \label{eq:TaTbcocycle}
    T_\alpha^{-1} T_\beta T_\beta^{-1} 
    T_\gamma T_\gamma^{-1} T_\alpha = \id
\end{equation}
whence 
$\Exp(t_{\alpha\beta}) \star \Exp(t_{\beta\gamma}) \star
\Exp(t_{\gamma\alpha})$ 
has to be a central element and thus \emph{constant} on the triple
overlap. Since for the star exponential one has a BCH formula denoted
by $\Exp(f)\star\Exp(g) = \Exp(f \circ_\star g)$ we find that
\begin{equation}
    \label{eq:tabc}
    t_{\alpha\beta\gamma} = 
    t_{\alpha\beta} \circ_\star
    t_{\beta\gamma} \circ_\star
    t_{\gamma\alpha} \in \mathbb{C}[[\lambda]]
\end{equation}
has to be \emph{constant} on 
$\mathcal{O}_\alpha \cap \mathcal{O}_\beta \cap \mathcal{O}_\gamma$ as
well. Thus we have a \v{C}ech cochain and it turns out that it is even
a cocycle. Thereby it defines a \v{C}ech cohomology class
$[t_{\alpha\beta\gamma}] \in \CH^2 (M, \mathbb{C})[[\lambda]]$. This
class coincides under the usual identification of 
$\CH^2(M, \mathbb{C})$ with $\HdR^2(M, \mathbb{C})$ with the relative
class, i.e.\  $[t_{\alpha\beta\gamma}] = t(\star', \star)$.

Now let us turn to (strong) Morita equivalence of star products. If
$\star$ and $\star'$ are Morita equivalent by some equivalence
bimodule, then it has to be a finitely generated and projective
$\star$-right module. Thus it is a deformed vector bundle 
$E \to M$ with deformation $\bullet$ and $\qh$. But since we want
$(C^\infty(M)[[\lambda]],\star')$ to be $^*$-isomorphic to
$\End(\Gamma^\infty(E)[[\lambda]], \bullet, \qh)$
it is necessary that $E$ is a \emph{line bundle} $L \to M$ since only
in this case $\Gamma^\infty (\End(L)) \cong C^\infty(M)$. Thus let 
$L \to M$ be a line bundle with Hermitian fiber metric and
corresponding deformations quantization $\bullet$ and $\qh$ with
respect to $\star$ such that $\star'$ is the star product for
$M$ such that
\begin{equation}
    \label{eq:MEfirststep}
    (C^\infty(M)[[\lambda]], \star') \cong
    \End(\Gamma^\infty(L)[[\lambda]], \bullet, \qh).
\end{equation}
We already know that it does not matter which deformation $\bullet$
and $\qh$ we actually choose since all are equivalent anyway. 
Moreover, we can choose within the \emph{isomorphism} class of
$\star'$ a representative such that we have a bimodule
deformation. This corresponds to the distinguished \emph{equivalence}
class as discussed in Theorem~\ref{theorem:defprojmod}. 
The subtle point here is the fact that an isomorphism of star products
starts in zeroth order with the pull-back by some diffeomorphism which
is not local. At the moment we want to exclude this non-locality and
consider this additional freedom later. For such $\star'$ we know from
Section~\ref{sec:ltog} that we can find \emph{local} operators
\begin{equation}
    \label{eq:Talpha}
    T_\alpha = \id + \sum_{r=1}^\infty \lambda^r T_\alpha^{(r)}: 
    C^\infty(\mathcal{O}_\alpha)[[\lambda]] \to
    C^\infty(\mathcal{O}_\alpha)[[\lambda]]
\end{equation}
such that $T_\alpha(f)$ is the local matrix of a global
$\bullet$-right linear endomorphism $T(f)$ of the deformed line
bundle. Moreover, $\star'$ can be obtained by 
$T(f \star' g) = T(f)\circ T(g)$. On the other hand, from
\eqref{eq:niceT} we know
\begin{equation}
    \label{eq:TalphaNice}
    T_\alpha (f \star' g) = 
    T_\alpha(f) \star T_\alpha (g).
\end{equation}
Up to now we have not yet used that $\star$ was a symplectic star
product. Thus we have the following statement, which is non-trivial in
the Poisson case:
\begin{lemma}
    \label{lemma:localequiv}
    Let $\star$ and $\star'$ be Morita equivalent star products on a
    Poisson mani\-fold with a deformed equivalence bimodule. Then
    $\star$ and $\star'$ are locally equivalent. In particular, they
    deform the same Poisson bracket, see also
    \cite[Lem.~3.4]{bursztyn:2001:pre} 
\end{lemma}
Since $T(f)$ is a global $\bullet$-right linear endomorphism we
conclude with \eqref{eq:endooverlap} that
\begin{equation}
    \label{eq:TaTbphiab}
    T_\alpha^{-1} T_\beta (f) 
    = \qphi_{\alpha\beta} \star f \star \qphi_{\beta\alpha}.
\end{equation}
Since in zeroth order $\qphi_{\alpha\beta}$ is just
$\phi_{\alpha\beta}$ we can write them as star exponentials 
$\qphi_{\alpha\beta} = \Exp(t_{\alpha\beta})$ with some
$t_{\alpha\beta} = \sum_{r=0}^\infty \lambda^r t_{\alpha\beta}^{(r)}$,
where $\eu^{t_{\alpha\beta}^{(0)}} = \phi_{\alpha\beta}$. Thus in the
symplectic case the relative class of $\star'$ and $\star$ is given by
the class $t(\star',\star) = [t_{\alpha\beta\gamma}]$ coming form
these $t_{\alpha\beta}$ as in \eqref{eq:tabc}. On the other hand, this
class can be computed explicitly. Since the $\qphi_{\alpha\beta}$
satisfy the cocycle identity \eqref{eq:cocycle} we have
\begin{equation}
    \label{eq:tabcintegral}
    t_{\alpha\beta\gamma} = 
    t_{\alpha\beta} \circ_\star
    t_{\beta\gamma} \circ_\star
    t_{\gamma\alpha} = 2 \pi \im \, n_{\alpha\beta\gamma}
    \quad
    \textrm{with}
    \quad
    n_{\alpha\beta\gamma} \in \mathbb{Z}.
\end{equation}
In particular, $t_{\alpha\beta\gamma}$ does \emph{not} depend on
$\lambda$ at all. But this shows that it coincides with its classical
limit which is given by 
$t_{\alpha\beta}^{(0)} + t_{\beta\gamma}^{(0)} +
t_{\gamma\alpha}^{(0)}$. Since the $t_{\alpha\beta}^{(0)}$ where just
the logarithms of the transitions functions $\phi_{\alpha\beta}$ of
the undeformed line bundle we conclude that 
$t(\star',\star) = 2\pi\im \, c_1(L)$ yields the \emph{Chern class} of
$L$. On the other hand, if $t(\star',\star)$ is a
\emph{$2\pi\im$-integral class} then there exists a line bundle with
fiber metric which we can deform. Then
Theorem~\ref{theorem:defprojmod} together with
Lemma~\ref{lemma:deffull} and Theorem~\ref{theorem:StrongylFull}
shows that $\star$ and $\star'$ are strongly Morita equivalent in this
case. Thus we have the following characterization  
\cite[Thm.~3.1 and Cor.~3.2]{bursztyn.waldmann:2001:pre}.
\begin{theorem}
    \label{theorem:mestars}
    Let $(M, \omega)$ be a symplectic manifold with Hermitian star
    products $\star$ and $\star'$.
    \begin{enumerate}
    \item $\star'$ is $^*$-equivalent to the star product induced by
        the deformation quantization of a line bundle $L \to M$ with
        respect to $\star$ if and only if
        \begin{equation}
            \label{eq:tceins}
            t(\star',\star) = 2 \pi \im \, c_1 (L).
        \end{equation}
    \item Two star products $\star$ and $\tilde{\star}$ are strongly
        Morita equivalent if and only if there exists a diffeomorphic
        star product $\star'$ to $\tilde{\star}$ such that
        \begin{equation}
            \label{eq:mestars}
            t(\star', \star) \in 2 \pi \im \, \HdR^2 (M, \mathbb{Z}).
        \end{equation}
   \end{enumerate}
\end{theorem}
Here $\HdR^2(M, \mathbb{Z})$ denotes the image of the integer
cohomology in the complex-valued deRham cohomology under the usual
inclusion map. Note also that we had to re-implement the freedom of a
diffeomorphism.
\begin{remark}
    There exists also a characterization of Morita equivalent star
    products in the Poisson case. If one formally `inverts' the
    condition 
    $\lambda c(\star') = \lambda c(\star) + 2\pi\im \lambda c_1(L)$
    one obtains a geometric series for the relation between the
    corresponding Poisson structures $\pi$ and $\pi'$. Here one has to
    observe that $c(\star)$ starts in order $\lambda^{-1}$ with the
    symplectic form $\omega$. Thus one obtains a relation between the
    formal deformations of the original Poisson bracket in
    Kontsevich's sense \cite{kontsevich:1997:pre}, which are known to
    classify star products in the Poisson case. That this heuristic
    reasoning is actually correct is a result of Jur\v{c}o, Schupp,
    and Wess \cite{jurco.schupp.wess:2001b:pre}.
\end{remark}

In the above construction the inner products only played a minor
role. Indeed, we can also consider the Morita equivalence of star
products from a purely ring-theoretical point of view. It turns out 
that for Hermitian star products this notion of Morita equivalence
\emph{coincides} with the strong Morita equivalence as $^*$-algebras
over $\mathbb{C}[[\lambda]]$, see also the discussion in
\cite{bursztyn.waldmann:2001:pre}. We have emphasized here the aspect
of strong Morita equivalence since this implies the equivalence of the
categories of $^*$-representations on pre-Hilbert spaces, which
from a physical point of view is the most interesting consequence of
Morita equivalence. Nevertheless, there are also several other
consequences.

Recall that the \emph{Picard group} $\Pic(\mathcal{A})$ of an algebra
is the group of equivalence classes of self-Morita equivalences. In
the classical case we have 
$\Pic(C^\infty(M)) = \CH^2(M, \mathbb{Z}) = \Pic(M)$, where the
(geometric) Picard group of a manifold is the group of equivalence
classes of line bundles with group structure from the tensor
product. From \cite{bursztyn:2001:pre} one knows that the Picard group
$\Pic(\mathcal{A})$ acts canonically on the moduli space of
equivalence classes of deformations of $\mathcal{A}$. Thus
Theorem~\ref{theorem:mestars} can also be understood as an explicit
computation of this group action: For the star products on a
symplectic manifold the Picard group acts on the equivalence classes
just by adding the $2 \pi \im$-multiple of the Chern class of the
corresponding line bundle in the zeroth order of $\lambda$ to the
characteristic class. 
\begin{corollary}
    \label{corollary:moritaI}
    Let $(M, \omega)$ be a symplectic manifold.
    Then action $\Phi$ of the Picard group $\Pic(M)$ on the
    equivalence classes of star products $\Def(M, \omega)$ is
    given by
    \begin{equation}
        \label{eq:PicAction}
        \Phi: \Pic(M) \times \Def(M, \omega) 
        \ni ([L], c(\star)) \; \mapsto \; 
        c(\star) + 2\pi\im\, c_1(L) \in \Def(M,\omega),
    \end{equation}
    with the identification~\eqref{eq:charclass}.
\end{corollary}
In the case where the characteristic class of $\star$ has a vanishing
zeroth order term we can also compute the Picard group of $\star$,
since in this case the additional diffeomorphism in \eqref{eq:mestars}
does not spoil the argument.
\begin{corollary}
    \label{corollary:moritaII}
    Let $(M, \omega)$ be a symplectic manifold with star product
    $\star$ such that $c(\star)$ has vanishing zeroth order term.
    \begin{enumerate}
    \item The Picard group of $(C^\infty(M)[[\lambda]], \star)$ is
        given by the group of equivalence classes of flat line bundles
        over $M$. 
    \item $\Gamma^\infty(L)$ admits a $\star$-bimodule structure
        deforming the classical one if and only if $L$ is flat, 
        i.e.\  $c_1(L) = 0$.
    \end{enumerate}
\end{corollary}
\begin{remark}
    \label{remark:dirac}
    A very explicit physical interpretation of the above
    characterization of Morita equivalent star products can be achieved
    in the case of star products on cotangent bundles $T^*Q$. Here one
    explicitly knows various star products, in particular a star
    product $\star^B$ as quantization of a charged particle in
    presence of a magnetic field $B$ on the configuration space
    $Q$. As usual $B$ is viewed as a closed two-form. 
    One knows that $c(\star^B)$ equals the class of $\im B$). Then
    $\star^B$ is Morita equivalent to $\star$ without $B$ 
    if and only if the magnetic field satisfies
    \emph{Dirac's quantization/integrality condition for the magnetic
      monopole charges} 
    described by $B$. Thus Dirac's condition turns out to
    imply \emph{equivalence} of the whole representation theory and
    not just the existence of a particular `Schr\"odinger-like'
    representation on section of the line bundle defined by $B$. In
    fact, the latter can be understood as Rieffel induced by the
    ordinary Schr\"odinger representation. Thus the above developed
    framework proves to give physically interesting and non-trivial
    applications, see \cite{bursztyn.waldmann:2001:pre} for a further
    discussion.
\end{remark}

%
% the bibliography
%
%
%\begin{footnotesize}
%    \bibliographystyle{ewde}
%    \bibliography{articles,books,preprints,misc}
%\end{footnotesize}

\end{document}